\documentclass[11pt]{article}
\usepackage{amssymb,amsmath,color}
\usepackage[dvips]{epsfig}
\usepackage{graphicx}
\usepackage[numbers,sort&compress]{natbib}
\usepackage{hyperref}

\usepackage{changes} 

\definechangesauthor[color=blue]{JR}

\definecolor{colorJRblue}{rgb}{0.,0.,1.}
\definecolor{colorJRred}{rgb}{1.,0.,0.}

\newtheorem{Lemma}{Lemma}
\newtheorem{Theorem}{Theorem}

\newtheorem{Corollary}{Corollary}
\newtheorem{Remark}{Remark}

\newenvironment{Proof}[1][.]%
 {\begin{trivlist}\item[]\textbf{Proof#1 }}%
 {\hspace*{\fill}$\rule{0.3\baselineskip}{0.35\baselineskip}$\end{trivlist}}

\def\eps{\varepsilon}
\def\rme{\mathrm{e}}
\def\rmi{\mathrm{i}}
\def\rmt{\mathrm{t}}
\def\rmd{\mathrm{d}}
\def\R{\mathbb{R}}
\def\C{\mathbb{C}}
\def\calL{\mathcal{L}}
\def\calF{\mathcal{F}}
\def\tr{\mathrm{tr}}

\def\cph{c_\mathrm{ph}}
\def\calO{\mathcal{O}}
\def\sgn{\mathrm{sgn}}
\def\m{\boldsymbol{m}}
\def\3{\boldsymbol{\hat{e}}_3}
\def\tp{\tilde p}
\def\op{\overline{p}}
\def\oq{\overline{q}}

\def\cc{c_{\rm cp}}
\def\omk{\Gamma}
\def\tGamma{\tilde\Gamma}
\def\tgamma{\tilde\gamma}
\def\tbeta{\tilde\beta}
\def\ttbeta{\bar{\beta}}
\def\tomega{\tilde\omega}
\def\tOmega{\varOmega}
\def\tmu{\tilde\mu}
\def\tt{{\tilde t}}

\def\mz{m_{3}^0}

\def\Im{\mathrm{Im}}
\def\Re{\mathrm{Re}}

\begin{document}
\begin{center}
{\huge Pattern formation in axially symmetric Landau-Lifshitz-Gilbert-Slonczewski equations}

\today

C. Melcher\footnote{Lehrstuhl I f\"ur Mathematik and JARA-FIT, RWTH Aachen University, 52056 Aachen, Germany, melcher@rwth-aachen.de}  \& J.D.M. Rademacher\footnote{Fachbereich 3 -- Mathematik, Universit\"at Bremen, Postfach 33 04 40, 28359 Bremen, Germany, rademach@math.uni-bremen.de}
 
\end{center}

\begin{abstract}
The Landau-Lifshitz-Gilbert-Slonczewski equation describes magnetization dynamics in the presence
of an applied field and a spin polarized current. In the case of axial symmetry and with focus on one space dimension, we investigate the emergence of space-time patterns in the form of wavetrains and coherent structures, whose local wavenumber varies in space. A major part of this study concerns existence and stability of wavetrains and of front- and domain wall-type coherent structures whose profiles asymptote to wavetrains or the constant up-/down-magnetizations. For certain polarization the Slonczewski term can be removed which allows for a more complete charaterization, including soliton-type solutions. Decisive for the solution structure is the polarization parameter as well as size of anisotropy compared with the difference of field intensity and current intensity normalized by the damping. 
\end{abstract}

\section{Introduction} \label{s:intro}

%
%
%
%
%

This paper concerns the analysis of spatio-temporal pattern formation for the axially symmetric Landau-Lifshitz-Gilbert-Slonczewski equation for which the applied magnetic field and current are aligned with or orthogonal to the material anisotropy. In one space dimension we thus consider
\begin{equation}\label{e:llg}
\partial_t \m = \m\times \Big[ \alpha \partial_t \m - \partial_x^2 \m + (\mu m_3 - h)\3  + \frac{ \beta}{1+ \cc \; m_3}\m \times  \3 \Big]
\end{equation}
as a model for the magnetization dynamics $\m=\m(x,t) \in \mathbb{S}^2$ (i.e. $\m$ is a direction field) driven by an external field 
$\boldsymbol{h}=h\3$ and current $\boldsymbol{j}=\frac{ \beta}{1+ \cc \; m_3}\3$ with polarization parameter $\cc \in (-1,1)$. The parameters $\alpha > 0$ and $\mu \in \R$
are the Gilbert damping factor and the anisotropy constant, respectively. A brief overview of the physical background and interpretation of terms is given below in Section \ref{s:model}.  

\begin{figure}
\begin{center}
\begin{tabular}{cc}
\includegraphics[width=0.2\textwidth]{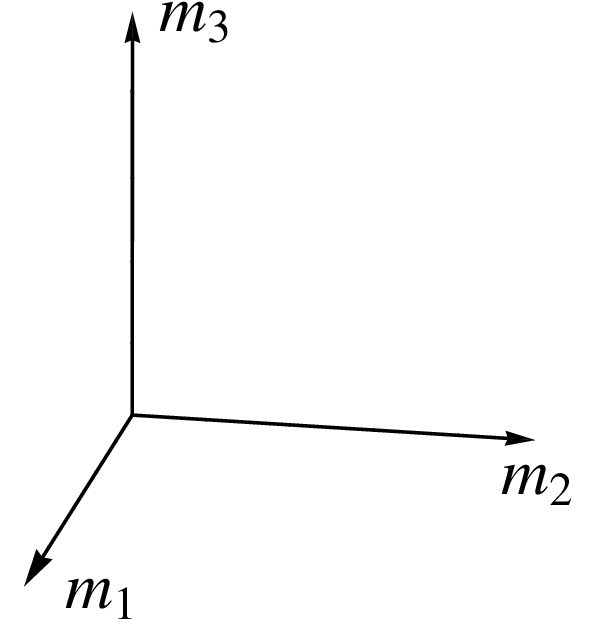}
\includegraphics[width=0.25\textwidth]{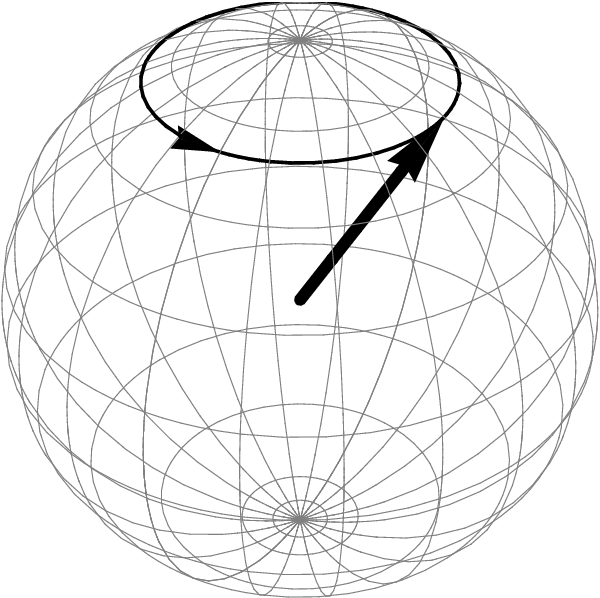}
& \hspace{5mm} \includegraphics[width=0.35\textwidth]{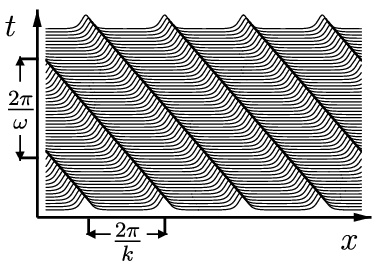}\\
(a) & (b)
\end{tabular}
\caption{Illustration of a wavetrain profile $\m(\varphi)$ (a) in the 2-sphere showing their constant altitude and (b) as a space-time plot of, e.g., $m_2$. In (a) the thick arrow represents $\m=(m,m_3)$, the thin line its trajectory as a function of $\varphi=kx-\omega t$.}
\label{f:wavetrain}
\end{center}
\end{figure}

The constant up- or down-magnetization states $\m=\pm \3$ are always steady states of \eqref{e:llg} and magnetic domain walls spatially separating these states are of major interest. While the combination of field and current excitations gives rise to a variety of pattern formation phenomena, see e.g. \cite{Bertotti:10, Hoefer:05, Hoefer:13, Kravchuk:13, LCC:14}, not much mathematically rigorous work is available so far, in particular for the dissipative case $\alpha>0$ that we consider. The case of axial symmetry is not only particularly convenient from a technical perspective. It offers at the same time valuable insight in the emergence of space-time patterns and displays strong similarities to better studied dynamical systems such as real and complex Ginzburg-Landau equations. In this framework we examine the existence and stability of wavetrain solutions of \eqref{e:llg}, i.e., 
solutions of the form
\[
\m(x,t) = e^{\rmi (kx-\omega t)} \m_0, 
\] 
where the complex exponential acts on $\m_0 \in \mathbb{S}^2$ by rotation about the $\3$-axis (cf. Figure~\ref{f:wavetrain} for an illustration). In the special case $\cc=0$ it turns out that \eqref{e:llg} can be transformed to the variational LLG-equation with $\beta=0$: in a rotating frame about the $m_3$-axis with frequency $-\beta/\alpha$ the current dependent term vanishes and $h$ changes to $h-\beta/\alpha$; see \S\ref{s:symc0}. This allows for a complete characterization of wavetrains and their $L^2$-stability. 

\begin{figure}
\begin{center}
\begin{tabular}{cc}
\includegraphics[width=0.4\textwidth]{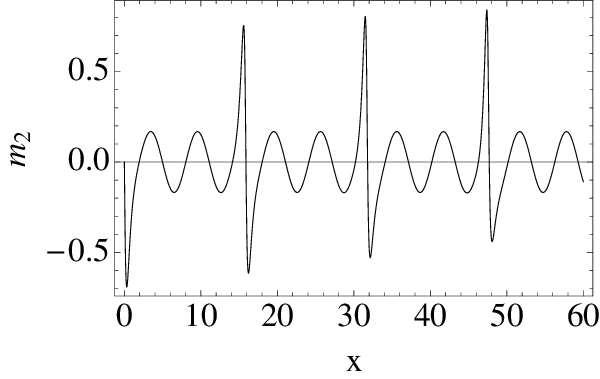} &
\includegraphics[width=0.4\textwidth]{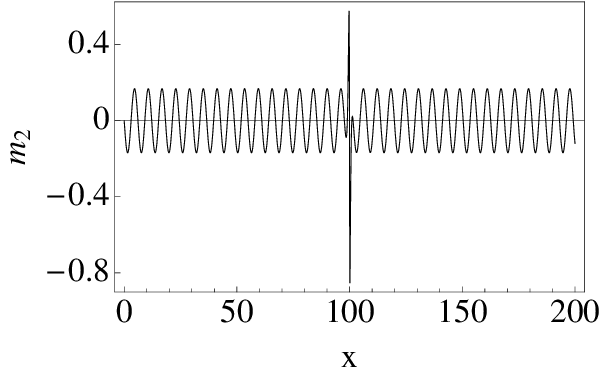}\\
\includegraphics[width=0.4\textwidth]{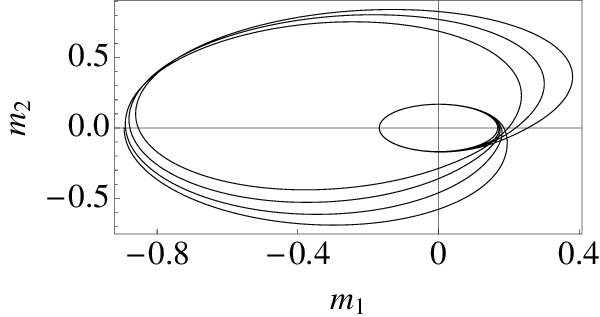} &
\includegraphics[width=0.4\textwidth]{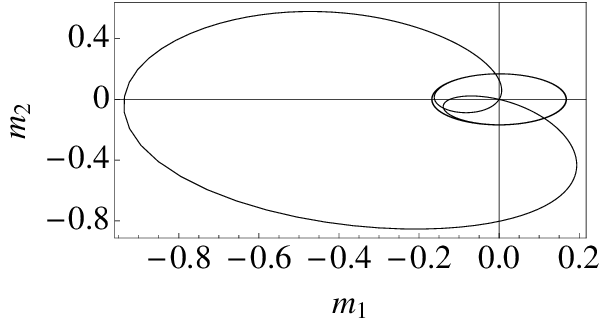}\\
(a) & (b)
\end{tabular}
\caption{Plots of coherent structure profiles for supercritical anisotropy with $\mu=7$, $h-\Omega=-1, \Omega=\beta/\alpha$, $\cc=0$, and first integral $C=1$, cf.~\eqref{e:coh-c0-q}. (a) A quasi-periodic solution that maps to a solution with period $\approx 16$ in the reduced equations \eqref{e:ode-c0-red}.  (b) The same solution type with period $200$ closer to a soliton-type solution with wavetrain as its asymptotic state. }
\label{f:cohex}
\end{center}
\end{figure}

\begin{figure}
\begin{center}
\begin{tabular}{cc}
\includegraphics[height=50mm]{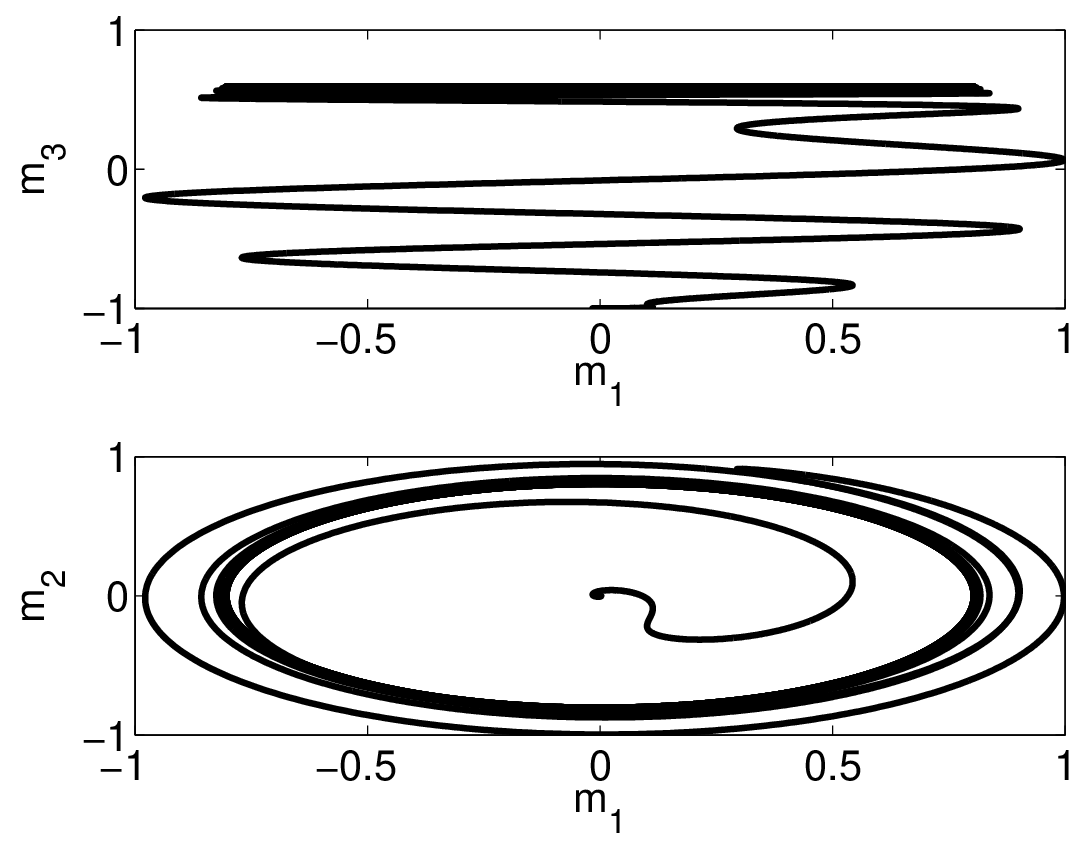} & 
\includegraphics[height=50mm]{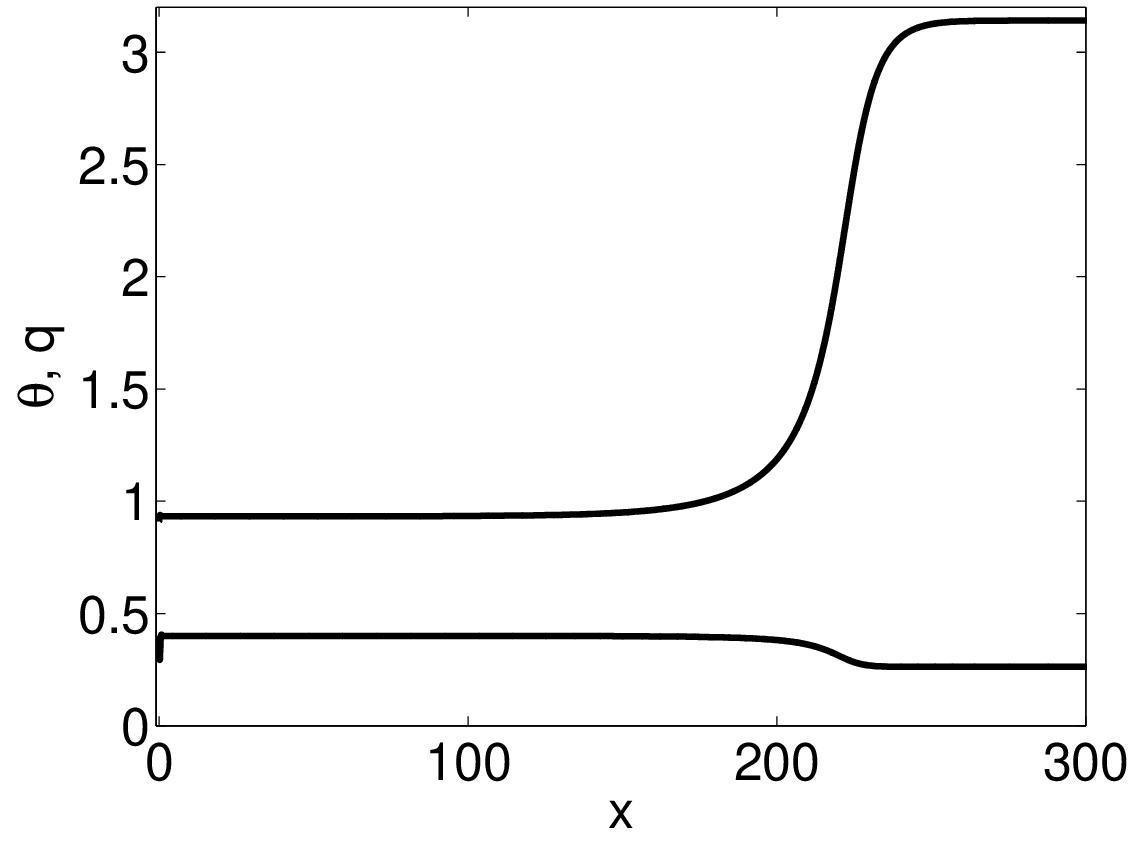}\\
(a) & (b)
\end{tabular}
\caption{Profile of a `fast' front connecting the wavetrain and the unstable $-\3$ computed with the coherent structure ODE guided by the asymptotic prediction of equation \eqref{e:superslow}. Here $\mu=1, h=0.5, s=5, \Omega=2, \cc=0$ and the asymptotic wavetrain on the left has wavenumber $k=\Omega/s=0.4$, and is spectrally stable.}\label{f:fast}
\end{center}
\end{figure}

\begin{figure}
\begin{center}
\begin{tabular}{ccc}
\includegraphics[width=0.3\textwidth]{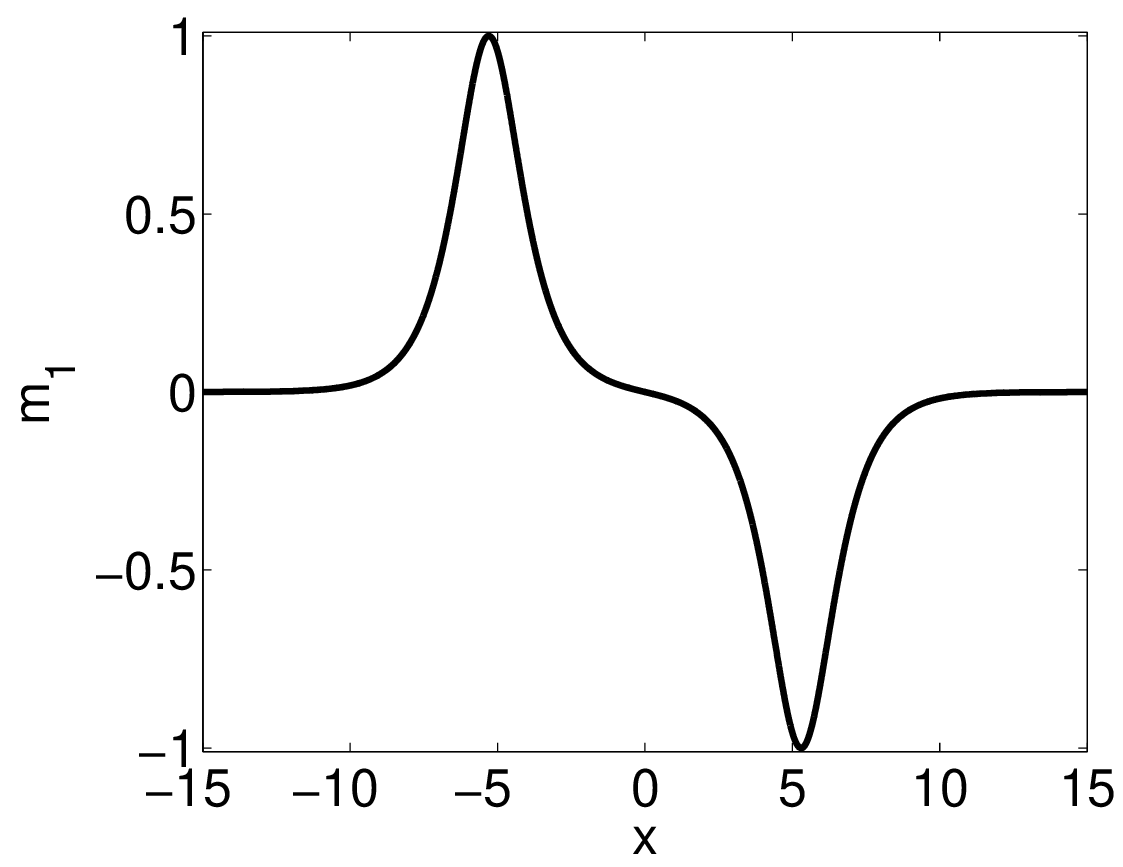}&
\includegraphics[width=0.3\textwidth]{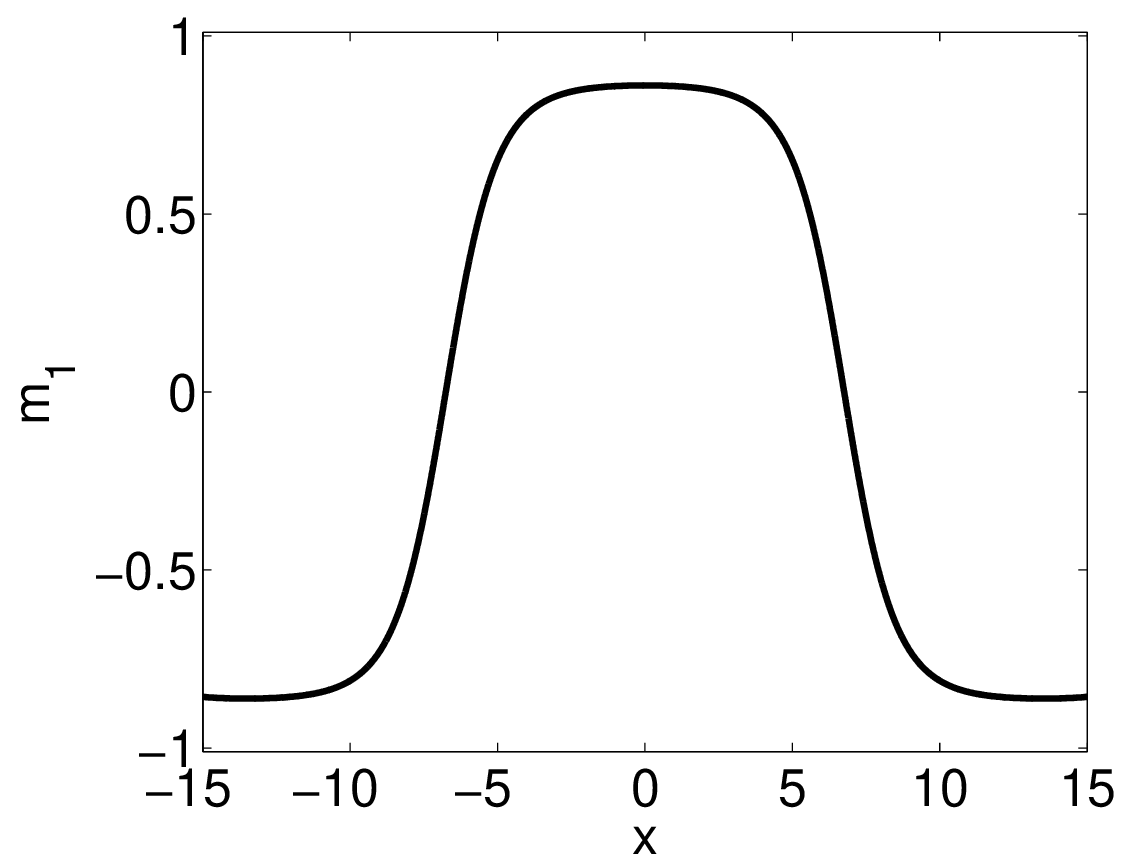}&
\includegraphics[width=0.3\textwidth]{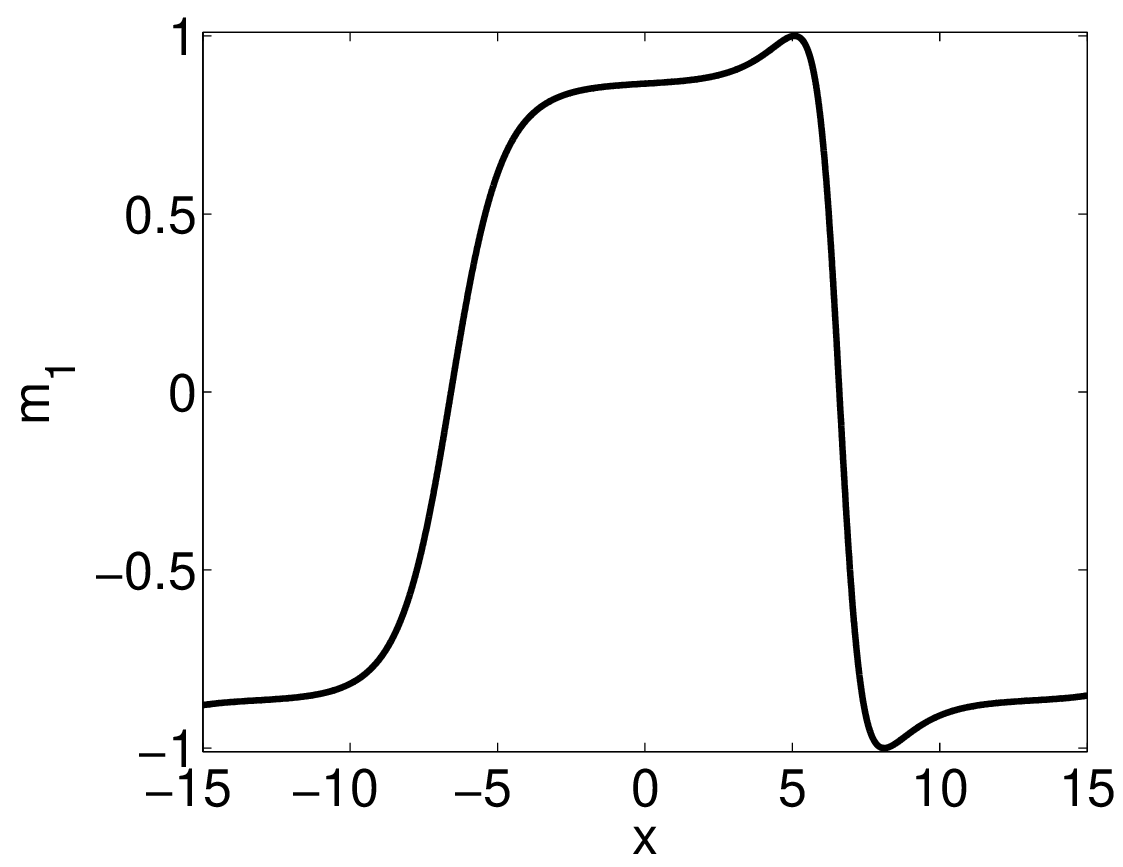}\\
\includegraphics[width=0.3\textwidth]{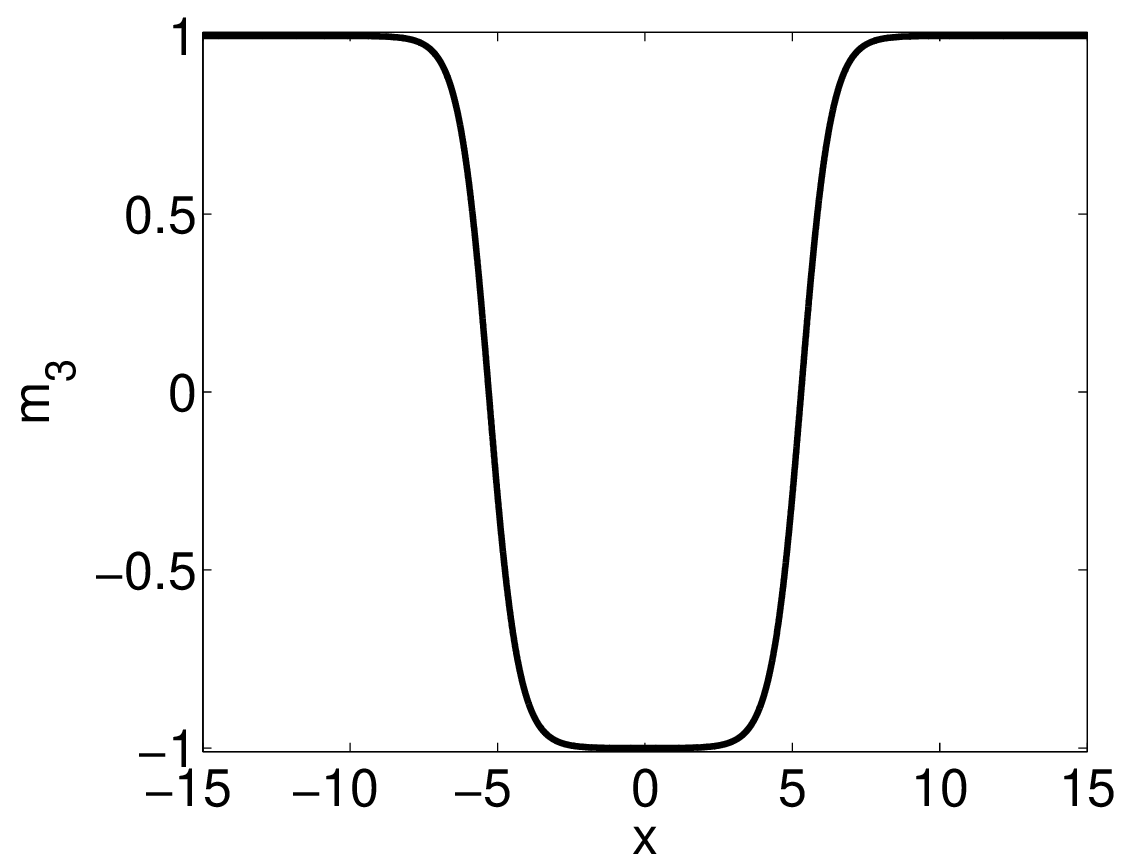}&
\includegraphics[width=0.3\textwidth]{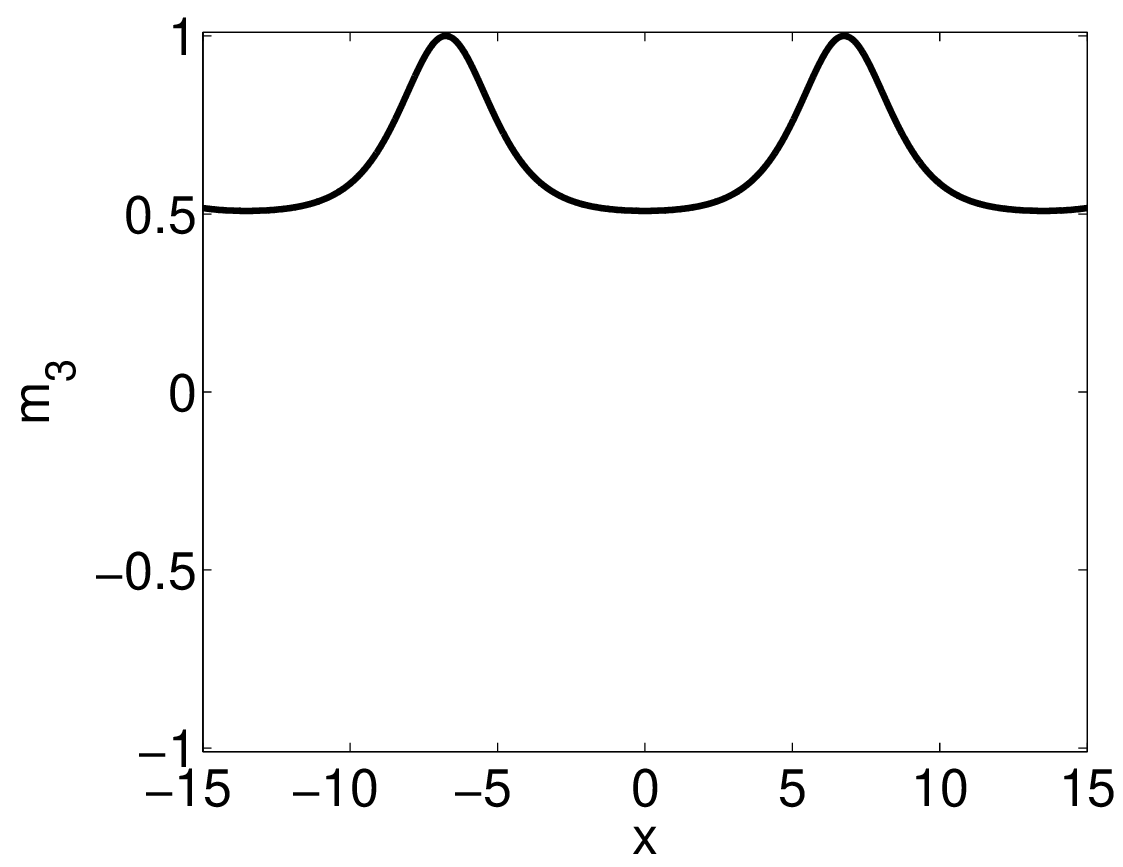}&
\includegraphics[width=0.3\textwidth]{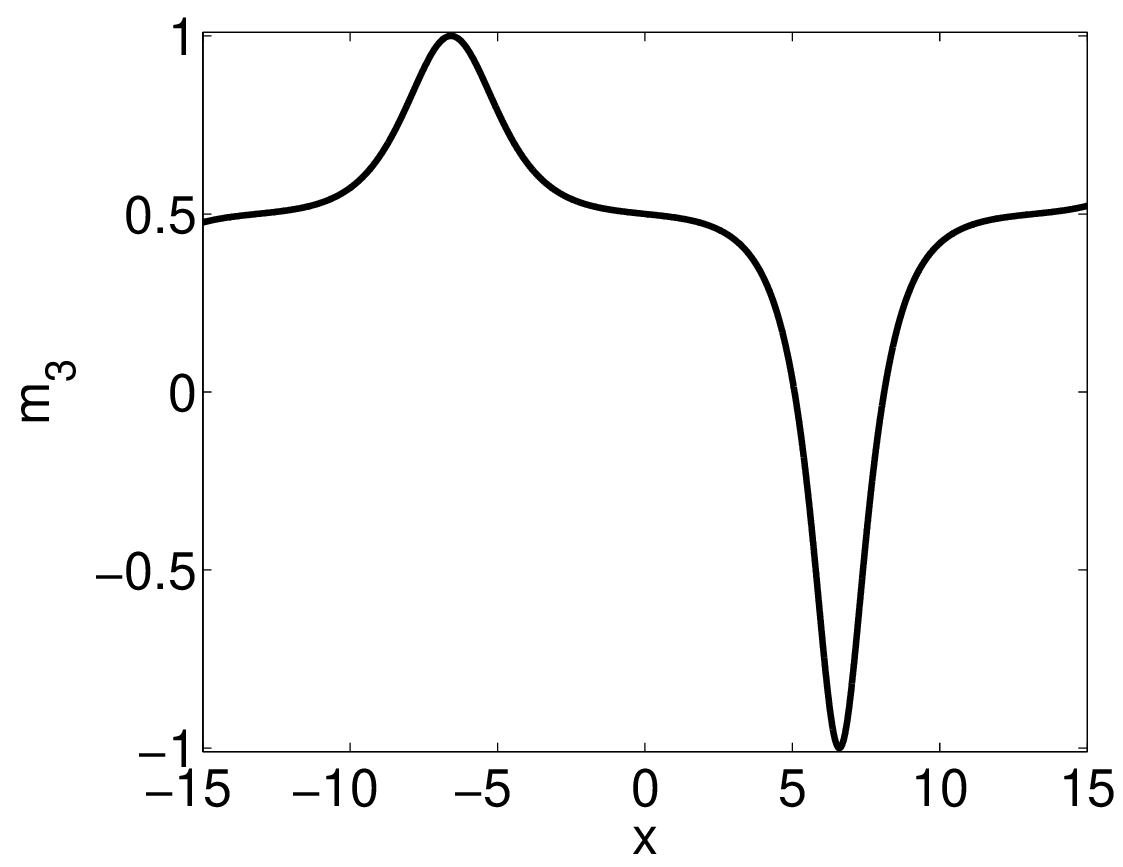}\\
(a) & (b) & (c)
\end{tabular}
\caption{Snapshots of sample homogeneous coherent structures with spatially periodic profiles, having $\rmd \phi/\rmd \xi = q=0$. Compare Figure~\ref{f:phaseplane-c0-q0}. Here $\cc=0,\alpha=1$ and $\Omega=\beta$ is arbitrary. 
(a) Near a pair of domains walls ($\mu=-1,h=10^{-4}-\beta$). (b) Near an upward `phase slip soliton' with plateaus near the oscillation at $m_3=h-\beta$ ($\mu=1, h=0.8-\beta$), and (c) near an upward-downward pair of such solitons ($\mu=1,h=0.8-\beta$). }
\label{f:profiles-c0-q0}\end{center}
\end{figure}

We also investigate the existence of coherent structure solutions which are locally in space of wavetrain form
\begin{equation}\label{e:coherent}
\m(x,t) = e^{i \varphi(x,t)} \m_0(x-st) \quad \text{where} \quad  \varphi(x,t)= \phi(x- s t) + \Omega t
\end{equation}
such that $\m_0(\xi) = (\sin \theta(\xi) ,0, \cos \theta(\xi) )$. Samples are plotted in Figures~\ref{f:cohex}, \ref{f:fast}, \ref{f:profiles-c0-q0}. In the variational case $\cc=0$ we completely characterize the existence of small amplitude coherent structures and stationary ($s=0$) ones, which in fact correspond to standing waves in the above rotating frame. Through the coherent structure viewpoint we recover a family of `homogenous' domain wall type solutions of arbitrary velocity, having no azimuthal profile, i.e., constant $\phi$ and thus vanishing local wavenumber $\rmd \varphi/\rmd x$. For general $\cc$ and large speeds, $|s|\gg 1$, we prove existence of a family of more general front-type coherent structures with nontrivial local wavenumbers, which can also form a spatial interface between $\pm\3$ and wavetrains. The analysis of these kinds of solutions is inspired by and bears similarities with that of the real and complex Ginzburg-Landau equations.

More specifically, the parameter space for existence of wavetrains and coherent structures is largely organized by the stability of the equilibrium states $\m=\pm\3$. The nature of bifurcations that we find motivates the following notions to organize the parameter space of \eqref{e:llg}: We refer to parameters as being
\begin{itemize}
\item `supercritical' if $\pm\3$ are both unstable 
\item `subcritical' if $\pm\3$ have different stability 
\item `subsubcritical' if $\pm\3$ are both stable.
\end{itemize}
From a physical viewpoint $\mu$ and $\alpha$ are material specific, while $h,\beta$ are control parameters. In our exposition we choose $\mu$ as a primary parameter and speak of super-, sub- or subsubcritical anisotropy; one may also choose $\beta$ or $h$ at the price of less convenient conditions.

Our results may be  summarized somewhat informally as follows.

\paragraph{The up- and down-magnetization equilibria $\pm\3$.} (Lemma~\ref{l:conststab})
Let $\beta^\pm:=\beta/(1\pm\cc)$. The constant state $\m=\3$ is strictly stable if and only if $\mu < h-{\beta^+}/\alpha$ and $\m=-\3$ if and only if $\mu <-(h-{\beta^-}/\alpha)$. Instabilities are of Hopf-type for the essential spectrum with onset via spatially homogenous modes of frequency $\beta/\alpha$.
In other words, stability of $\pm\3$ changes when the difference of signed anisotropy $\pm\mu$ and the force balance $h-\frac{\beta^\pm}{\alpha}$, of magnetic field strength minus the ratio of current-polarization intensity and damping factor, changes sign. This corresponds to the well known instability threshold in the more broadly studied ODE for solutions that are homogeneous in space. Notably, for $\cc=0$ the anisotropy is subsubcritical precisely for $-\mu>|h-\beta/\alpha|$ (`easy-axis'), and supercritical precisely for $\mu>|h-\beta/\alpha|$ (`easy-plane'). 


\paragraph{Fast and small amplitude coherent structure.} (Theorems~\ref{t:coh-fast}, \ref{t:coh-small} and corollaries)
For each sufficiently large speed there exists a family of front-type coherent structures parametrized by the azimuthal frequency. Their profiles connect $\pm\3$ with each other or, if there are wavetrains, there are fronts connecting these and/or $\pm\3$ in the order of altitudes. An example is plotted in Figure~\ref{f:fast}. Small amplitude coherent structures are of front type and, for $\cc=0$, exist only for super- and subcritical anistropy.  

\paragraph{Wavetrains.} 
In the case $\cc=0$ (Theorems~\ref{t:wt-ex}, \ref{t:wt}) for each wavenumber $k\in\R$ at most one wavetrain exists, and moreover:
\begin{enumerate}
\item \emph{Supercritical anisotropy:} Wavetrains exist precisely for $k$ with $|k|>\sqrt{\mu+|h-\beta/\alpha|}$ or $0\leq |k| <\sqrt{\mu-|h-\beta/\alpha|}$. There is (explicitly known) $k_*\in(0,\sqrt{\mu-|h-\beta/\alpha|}$ such that all wavetrains with $|k|<k_*$ are stable and sideband unstable for $|k|>k_*$.
\item \emph{Subcritical anisotropy:} Wavetrain exists precisely for $k$ with $|k|>\sqrt{\mu+|h-\beta/\alpha|}$, but are all unstable.
\item \emph{Subsubcritical anisotropy:}
Wavetrains exist for all $k$,  but are all unstable.
\end{enumerate}

The overall picture for wavetrains of \eqref{e:llg} with $\cc=0$ can be viewed as a combination of those in a supercritical and a subcritical real Ginzburg-Landau equation; see Figure~\ref{f:wavetrainstab}.

\medskip
For general $\cc\in(-1,1)$ additional effects are (1) a nontrivial nonlinear dispersion relation $\omega(k)$ with nonzero group velocities $\frac{\rmd}{\rmd k}\omega(k)$, (2) the occurrence of `hyperbolic' and `elliptic'  bifurcation points of wavetrains and (3) coexistence of stable wavetrains and stable $\pm\3$. Wavetrains for $k^2>\mu$ are always unstable (Theorems~\ref{t:wt-ex}, \ref{t:wt-stab}), but for $\cc\neq 0$ wavetrains are potentially convectively but not absolutely unstable, though we do not investigate this here.

\medskip
\noindent\textbf{Domain walls for $\cc=0$.} (Theorem~\ref{t:domwall}) For any $\mu<0$ there exists a family of fronts whose spatial profiles connect $\pm\3$ with $\theta'=\sqrt{-\mu}\sin(\theta)$, and that are `homogeneous' in the sense that $q\equiv 0$ so there is no azimuthal profile. They corresponds to well known domains walls of the LLG-equation\footnote{After acceptance of the present manuscript for publication, we found these were also obtained in \cite{GRS:10}.} . Here we readily locate these within the coherent structure framework.

\medskip
\noindent\textbf{Stationary coherent structures for $\cc=0$.}  (Theorems~\ref{t:coh-c0-q0}, \ref{t:coh-c0}) 
\begin{enumerate}
\item \emph{Supercritical anisotropy:} For fixed parameters there exist various stationary coherent structures ($s=0$) including `homogeneous' ones (cf.\ Figure~\ref{f:profiles-c0-q0}). An interesting case of the latter is a symmetric pair of `phase slip' soliton-type coherent structures, whose spatially asymptotic states are the same spatially homogeneous oscillation ($k=0$), but the intermediate profile crosses either $\3$ or $-\3$, so that the asymptotic states differ azimuthally by $180^\circ$. There also exists a non-homogeneous soliton-type solution with asymptotic state being a wavetrain (cf.\ Figure~\ref{f:cohex}).
\item \emph{Sub- and subsubcritical anisotropy:} All stationary coherent structures have periodic profiles except a homogeneous phase slip soliton with spatially asymptotic state $\pm\3$ for $\sgn(h-\beta/\alpha)=\pm1$. 
\end{enumerate}


\paragraph{Higher space dimensions.}

The model for $N$ space dimensions has the second derivative with respect to $x$ in \eqref{e:llg} replaced by a Laplace operator $\sum_{j=1}^N \partial_{x_j}^2$. Wavetrain type solutions are then of the form
\[
\m(x,t) = \m_*(k \cdot x-\omega t), 
\]
where $k=(k_1,\ldots,k_N)$. Notably, for $k_j=0$, $2\leq j \leq N$ these are solutions from one space dimension extended trivially (constant) in the additional directions. 

Conveniently, the rotation symmetry (gauge invariance in the Ginzburg-Landau context), means that the analyses of $\pm\3$ and these wavetrains is already covered by that of the one-dimensional case: the linearization is space-independent and therefore there is no symmetry breaking due to different $k_j$.
Indeed, all relevant quantities are rotation symmetric, depending only on $k^2=\sum_{j=1}^N k_j^2$ or $\ell^2=\sum_{j=1}^N \ell_j^2$, where $\ell=(\ell_1,\ldots,\ell_N)$ is the Fourier wavenumber vector of the linearization. In particular, the instabilities occur simultaneously for all directions.

Concerning coherent structures, in higher space dimension the defining equation (see \eqref{e:ode} below) turns into an elliptic PDE in general.  The analysis in this paper only covers the trivial constant extension into higher dimensions.

\bigskip
This paper is organized as follows. In Section \ref{s:model}, the terms in the model equation \eqref{e:llg} and its well-posedness are discussed. Section~\ref{s:hopf} concerns the stability of the trivial steady states $\pm\3$ and in \S\ref{s:wt} existence and stability of wavetrains are analyzed. Section \ref{s:coh} is devoted to coherent structures.


\textbf{Acknowledgement.}  JR has been supported in part by the NDNS+ cluster of the Dutch Science Fund (NWO). We thank  the anonymous reviewers for suggestions that helped improve the manuscript, and Lars Siemer as well as Ivan Ovsyannikov for their critical reading.


\section{Review of Landau-Lifshitz-Gilbert-Slonczewski equations}\label{s:model}

The classical equation of dissipative magnetization dynamics, the Landau-Lifshitz-Gilbert equation \cite{Landau_Lifshitz:35, Gilbert:04} for unit vector fields $\m=\m(x,t) \in \mathbb{S}^2$, 
\[
\partial_t \m= \m \times \left( \alpha \, \partial_t \m - \gamma \,  \boldsymbol{h}_{\rm eff} \right).
\] 
features a damped precession of $\m$ around the effective field $\boldsymbol{h}_{\rm eff}= - \delta \mathcal E(\m)$, 
i.e., minus the variational derivative of the interaction energy $\mathcal E=\mathcal  E(\m)$. 
The gyromagnetic ratio $\gamma>0$ is a parameter which appears as the typical precession frequency. By rescaling time, 
one can always assume $\gamma=1$. The Gilbert damping factor $\alpha>0$ is a constant that can be interpreted dynamically as the inverse of the typical relaxation time.
It is useful to take into account that there are several equivalent forms of LLG.  Elementary algebraic manipulations taking into account that 
$
- \m \times \m \times \boldsymbol{\xi} = \boldsymbol{\xi} - (\m \cdot \boldsymbol{\xi}) \m
$
yield the so-called Landau-Lifshitz form
\begin{equation}\label{eq:LL_form}
(1+ \alpha^2) \partial_t {\m} =-  \m \times \left( \alpha \, \m \times    \boldsymbol{h}_{\rm eff}   + \boldsymbol{h}_{\rm eff}  \right),
\end{equation}
introduced in the original work \cite{Landau_Lifshitz:35}. 
In case $\alpha>0$, the energy $\mathcal{E}(\m)$ is not conserved but is a Lyapunov functional, i.e., more precisely 
(recall $\boldsymbol{h}_{\rm eff}= - \delta \mathcal{E}(\m)$)
\[
\frac{d}{dt} \mathcal{E}(\m(t))= - \alpha \|\partial_t \m(t)\|^2 \quad \text{or equivalently} \quad 
 \frac{d}{dt} \mathcal{E}(\m(t))= - \frac{\alpha}{1+\alpha^2} \| \m \times \boldsymbol{h}_{\rm eff} \|^2.
\]
Gilbert damping enables the magnetization to approach (spiral down to) a steady state, i.e. satisfying $\m \times \boldsymbol{h}_{\rm eff}=0$
(Browns equation), as $t \to \infty$. 

\paragraph{Spin-torque interaction.}
The system can be driven out of equilibrium conventionally by an external magnetic field $\boldsymbol{h}$ which appears as part of the effective field. 
In modern spintronic applications, magnetic systems are excited by spin polarized currents (with direction of polarization $\boldsymbol{\hat{e}}_p \in \mathbb{S}^2$) giving rise to a spin torque
\begin{equation}\label{eq:spin_torque}
 \m \times  \m \times {\bf j} \quad \text{where} \quad  \boldsymbol{j}= \beta \frac{\boldsymbol{\hat{e}}_p}{1+\cc  \; \m \cdot \boldsymbol{\hat{e}}_p}, 
\end{equation}
which has been introduced in  \cite{Berger:96, Slonczewski:96}. Here, the parameters 
$\beta>0$ and  $\cc \in (-1,1)$ depend on the intensity of the current and ratio of polarization \cite{Bertotti:08}. Typically we have $ \boldsymbol{\hat{e}}_p = \3$.
Accordingly, the modified Landau-Lifshitz-Gilbert equation, also called Landau-Lifshitz-Gilbert-Slonczewski equation (LLGS), reads
\begin{equation}\label{eq:LLG_ODE}
\partial_t \m = \m \times \left( \alpha \, \partial_t \m - \boldsymbol{h}_{\rm eff} + \m \times \boldsymbol{j}  \right).
\end{equation}
One may extend the notion of effective field to include current interaction by letting
 \[
 \boldsymbol{H}_{\rm eff} =\boldsymbol{h}_{\rm eff} - \m \times \boldsymbol{j},
 \] 
where the second term is usually called Slonczewski term. In this framework \eqref{eq:LLG_ODE} can also be written in the form
\eqref{eq:LL_form} with $\boldsymbol{h}_{\rm eff}$ replaced by $\boldsymbol{H}_{\rm eff}$.
Observe, however,  that the Slonczewski term (and hence $\boldsymbol{H}_{\rm eff}$) is in general non-variational and that the energy is no longer a Lyapunov functional. 
Introducing the potential 
 $\Psi(\m) = \frac{\beta}{\cc} \ln(1+\cc \, \m \cdot \boldsymbol{\hat{e}}_p)$ of $\boldsymbol{j}$ (for $\cc\neq 0$) reveals the
\textit{skew variational} structure
\[
 \m \times \left[ \alpha \, \partial_t \m+\delta \mathcal{E}(\m) \right] =-  \m \times \m \times \left[ \partial_t \m +\delta \Psi(\m)  \right],
\]
see \cite{BMS:08}.
%
In the micromagnetic model
the underlying interaction energies are integral functionals in $\m$ containing in particular exchange (Dirichlet) interaction, dipolar stray-field interaction, crystal anisotropy and Zeeman interaction with external magnetic field, see e.g. \cite{Hubert_Schaefer}. In this paper
we shall mainly focus on the spatially one-dimensional situation and consider energies of the form
\begin{equation} \label{eq:energy}
\mathcal{E}(\m) =  \frac{1}{2} \int   \left( |\partial_x \m|^2 + \mu  m_3^2\right) \, dx - \int \boldsymbol{h} \cdot \m \; dx.
\end{equation}
Here, $\boldsymbol{h} \in \R^3$ is a constant applied magnetic field. The parameter $\mu \in \R$ features \textit{easy plane} anisotropy for
$\mu>0$ and  \textit{easy axis} anisotropy for $\mu <0$, respectively, according to energetically preferred subspaces. This term
comprises crystalline and shape anisotropy effects. Shape anisotropy typically arises from stray-field interactions which prefer magnetizations tangential to the sample boundaries.  Hence $\mu >0$ corresponds to a thin-film perpendicular to the $\3$-axis whereas $\mu <0$ corresponds to a thin wire parallel to the $\3$-axis. 
The effective field corresponding to \eqref{eq:energy} reads
\begin{equation}
\boldsymbol{h}_{\rm eff}= \partial_x^2 \m - \mu m_3 \3 + \boldsymbol{h}.
\end{equation}
With the choices $\boldsymbol{h}= h\, \3$ and $\hat{e}_p= \3$, the Landau-Lifshitz-Gilbert-Slonczewski equation \eqref{eq:LLG_ODE} exhibits the aforemented rotation symmetry about the $\boldsymbol{\hat e}_3$-axis. The presence of a spin torque $\m \times \m \times \boldsymbol{j}$ exerted by a constant current may induce switching between magnetization states or magnetization oscillation 
\cite{Berkov_Miltat:08, Bertotti:05, Bertotti:08}. For the latter effect, the energy supply due to the electric current compensates the energy dissipation due to damping enabling a stable oscillation, called \textit{precessional states}.
In applications the typical frequency is in the range of GHz, so that a precessional state would basically act as a microwave generator. In the class of spatially homogeneous states,
precessional states are periodic orbits with $m_3=const.$ and of constant angular velocity $ \beta/\alpha$ when $\cc=0$. It is more subtle, however, to understand the occurrence and stability of spatially non-homogeneous precessional states. This is the theme of this paper.

\paragraph{Extensions and related work.} 

There is a wealth of literature studying the dynamics of related Landau-Lifschitz models with and without damping and axial symmetry and including effects other than spin-torque interaction and as general reference we mention the book \cite{BMS:08} as well as the review article \cite{Lakshmanan:11}. More specifically,  spatially non-trivial states and their stabililty have been considered in \cite{KNML:05}, where the spin-torque part of the effective field is replaced by a demagnetization term solving Maxwell's equation. Also coupled nano-oscillators of LLGS type have been considered widely, e.g. recently in \cite{Shaffer,SCL:15}. Recently, for a situation without axial symmetry, Turing patterns of spin states have been numerically found in \cite{LCC:14}. 

Non-symmetric variants of our equation  \eqref{e:llg} have been used e.g. in the description of the field  driven motion of a flat domain wall connecting antipodal steady states $m_3 = \pm 1$ as $x_1 \to \pm \infty$. A prototypical situation is the field driven motion of a flat Bloch wall in an uniaxial the bulk magnet governed by
\begin{equation} \label{eq:DWmotion}
\partial_t \m = \m\times \left(\alpha \partial_t \m - \partial_x^2 \m +  \mu_1 m_1 \boldsymbol{\hat e}_1 + ( \mu_3 m_3 \boldsymbol - h ) \boldsymbol{\hat e}_3 \right).
\end{equation}
In this case $\mu_1>0> \mu_3$, where $\mu_1$ corresponds to stray-field and $\mu_3$ to crystalline anisotropy.
Explicit traveling wave solutions were obtained in unpublished work by Walker,  see e.g. \cite{Hubert_Schaefer}, and reveal interesting effects such as the existence of a terminal velocity (called Walker velocity) and the notion of an effective wall mass. A mathematical account on Walker's explicit solutions and investigations on their stability, possible extensions to finite layers and curved walls can be found e.g. in \cite{Carbou:10, Melcher:04, Podio_Tomassetti:04}. Observe that our axially symmetric model is obtained in the limit $\mu_1 \searrow 0$. On the other hand, the singular limit $\mu_3 \to + \infty$ leads to trajectories confined to the $\{m_3=0\}$ plane (equator map), and can be interpreted as a thin-film limit. In suitable parameter regimes it can be shown that the limit equation is a dissipative wave equation governing the motion of N\'eel walls
 \cite{Capella:07, Melcher:03, Melcher:10}.

\paragraph{Well-posedness of LLGS.} 
It is well-known that Landau-Lifshitz-Gilbert equations and its variants have the structure of quasilinear parabolic systems. In the specific case of  \eqref{e:llg}, one has the
extended effective field $\boldsymbol{H}_{\rm eff} =\boldsymbol{h}_{\rm eff} - \m \times \boldsymbol{j}$, more precisely
\begin{equation}\label{eq:H_ext}
\boldsymbol{H}_{\rm eff} = \partial^2_x \m - f(\m)  \quad \text{where} \quad f(\m)= (\mu m_3 - h) \3 + \frac{\beta}{1+\cc m_3} \m \times \3.
\end{equation}
Hence the corresponding Landau-Lifshitz form \eqref{eq:LL_form} of  \eqref{e:llg} reads
\begin{equation}\label{eq:LL_form_1}
(1+\alpha^2) \partial_t \m =-  \m\times \Big[ \partial_x^2 \m -  f(\m) \Big]- \alpha \m \times  \m\times \Big[ \partial_x^2 \m -  f(\m) \Big].
\end{equation}
Taking into account
\begin{equation} \label{eq:div}
  \m\times \partial_x^2 \m =  \partial_x( \m\times \partial_x \m) \quad \text{and} \quad - \m \times  \m \times  \partial_x^2 \m =  \partial_x^2 \m + |\partial_x \m|^2 \m,
\end{equation}
valid for $\m$ sufficiently smooth and $|\m|=1$, one sees that \eqref{eq:LL_form_1} has the form
%
\begin{align}\label{e:quasi}
\partial_t \m = \partial_x \left( A(\m) \partial_x \m \right) + B(\m, \partial_x \m)
\end{align}
with analytic functions $A:\R^3 \to \R^{3 \times 3}$ and $B:\R^3 \times \R^{3} \to \R^3$ such that $A(\m)$ is uniformly elliptic for $\alpha>0$, in fact
\[  
\boldsymbol{\xi} \cdot A(\m)  \boldsymbol{\xi} = \frac{\alpha}{1+\alpha^2} | \boldsymbol{\xi}|^2
\quad \text{for all} \quad \boldsymbol{\xi} \in \R^3.
\]

Well-posedness results for $\alpha>0$ can now be deduced from techniques based on higher order energy estimates as in \cite{Melcher:11, Melcher_Ptashnyk} or maximal regularity and interpolation as in \cite{users_guide}. In particular, we shall rely on results concerning perturbations of wavetrains, traveling waves,
and steady states. Suppose $\m_\ast=\m_\ast(x,t)$ is a smooth solution of \eqref{e:llg} with bounded derivatives up to all high orders (only sufficiently many are needed) and 
$\m_0:\R \to \mathbb{S}^2$ is such that $\m_0- \m_\ast(\cdot,0) \in H^2(\R)$. Then there exist
$T>0$ and a smooth solution $\m: \R \times (0,T) \to  \mathbb{S}^2$ of \eqref{e:llg} such that
$\m-\m_\ast \in C^0([0,T); H^2(\R)) \cap C^1([0,T); L^2(\R))$ with
\[
\lim_{t \searrow 0}  \| \m(t)- \m_0\|_{H^2} =0 \quad \text{and} \quad  \lim_{t \nearrow T}  \| \m(t)- \m_\ast(\cdot,t)\|_{H^2} =\infty
\quad \text{if $T<\infty$}.
\]
The solution is unique in its class and the flow map depends smoothly on initial conditions and parameters.  

\medskip
Given the smoothness of solutions, we may compute pointwise $\partial_t |\m|^2 = 2\m \cdot \partial_t \m$, so that for $|\m|=1$ the cross product form of the right hand side of \eqref{eq:LL_form_1} gives $\partial_t |\m|^2=0$. Hence, the set of unit vector fields, $\{|\m|=1\}$, is an invariant manifold of \eqref{e:quasi} consisting of  the solutions to \eqref{e:llg} that we are interested in.

In addition to well-posedness, also stability and spectral theory for \eqref{e:quasi}, see, e.g., \cite{users_guide}, carry over to \eqref{e:llg}. In particular, the computations of $L^2$-spectra in the following sections are justified and yield nonlinear stability for strictly stable spectrum and nonlinear instability for the unstable (essential) spectrum.

\paragraph{Landau-Lifshitz-Gilbert-Slonczewski versus complex Ginzburg-Landau equations.}
Stereographic projection of \eqref{e:llg} yields 
\[
 (\alpha + {\rm i})  \zeta_t =  \partial_x^2 \zeta -
 \frac{2 \bar{\zeta} (\partial_x \zeta)^2}{1+|\zeta|^2} + \mu   \frac{(1- |\zeta|^2) \zeta }{1+ |\zeta|^2}  - (h + i \beta) \zeta
\quad \text{where} \quad 
\zeta = \frac{m_1+{\rm i} m_2}{1+m_3},
\]
valid for magnetizations avoiding the south pole.

Studying LLG-type equations via stereographic projection has a long history and has been employed in several of the aforementioned references, see, e.g., the review \cite{Lakshmanan:11} and the references therein.mThere is also a global connection between LLG and CGL in the spirit of the classical Hasimoto transformation \cite{Hasimoto:72}, which turns the (undamped) Landau-Lifshitz equation in one space dimension ($\boldsymbol{h}_{\rm eff}=\partial_x^2$) into the focussing cubic 
Schr\"odinger equation \cite{Lakshmanan:77,Zhakharov:79}. The idea is to disregard the customary coordinates representation and to introduce instead a pull-back frame on the tangent bundle along $\m$. In the case of $\mu=\beta=h=0$, i.e. $\boldsymbol{h}_{\rm eff}= \partial_x \m$, this leads to
\begin{equation}\label{eq:gauged_LLS}
 (\alpha + \rmi)   \mathcal{D}_t u   =    \mathcal{D}_{x}^2 u 
\end{equation}
where $u=u(x,t)$ is the complex coordinate of $\partial_x \m$ in the moving frame representation, and $\mathcal{D}_x$ and $\mathcal{D}_t$ are
covariant derivatives in space and time giving rise to cubic and quintic nonlinearities, see \cite{Melcher:11, Melcher_Ptashnyk} for details. 


\subsection{Symmetry and the variational structure for $\cc=0$}\label{s:symc0}

The aforementioned rotation symmetry of \eqref{e:llg} about the $m_3$-axis of all terms manifests as an equivariance of the right hand side of \eqref{e:quasi} with respect to any such rotation $R_\varphi$: let $\m=R_\varphi\widetilde\m$, then
\[
\partial_x \left( A(\m) \partial_x \m \right) + B(\m, \partial_x \m) = R_\varphi( \partial_x \left( A(\widetilde\m) \partial_x \widetilde\m \right) + B(\widetilde\m, \partial_x \widetilde\m)).
\]

 For a rotating frame $\m =e^{\rmi\tOmega t}\m$, write the rotation about the $\3$-axis as $R_{\tOmega t}$ and note that time derivatives become $\partial_t\m = R_{\tOmega t}\left(\tOmega R_{\tOmega t}^*R'_{\tOmega t}\widetilde\m + \partial_t\widetilde\m\right)$, where $R_{\tOmega t}^*R'_{\tOmega t}\m = \m\times\3$. Therefore, having $\cc=0$, \eqref{e:llg} is also an equation for $\widetilde\m$ with the parameter $\beta$ changed to $\beta+\alpha\tOmega$ (from $\partial_t \m$ within the brackets) and $h$ to $ h+\tOmega$ (from $\partial_t \m$ on the left hand side). In other words, for $\cc=0$, changing spin torque current has the same effect as changing the magnetic field.
 
Choosing $\tOmega=-\beta/\alpha$ yields $\beta=0$ in \eqref{e:llg}, which is therefore variational with respect to the energy \eqref{eq:energy} as discussed above. This has strong structural consequences for the coherent structures \eqref{e:coherent} and allows for a (largely) complete characterization. In particular, it turns out that the existence of coherent structures that are stationary ($s=0$), but not necessarily time independent, requires their superimposed azimuthal frequency $\Omega$ to satisfy $\Omega=\beta/\alpha$; see \S\ref{s:stat-coh}. The reduction to $\beta=0$ thus implies $\Omega=0$ and therefore turns the stationary coherent structures into standing waves and hence to time-independent solutions.

\section{Hopf instabilities of the steady states $\m=\pm \3$}\label{s:hopf}

As a starting point and to motivate the subsequent analysis of more complex patterns, let us consider the stability of the constant magnetizations $\pm\3$. It is well-known that a Hopf bifurcation of these states occurs in the ODE associated to \eqref{e:llg} in the absence of diffusion, that is, for spatially constant solutions. In the following, we account in addition for spatial dependence.

Use the shorthand $\beta^\pm=\beta/(1\pm \cc)$. Substituting $\m= \pm\3 +\delta \boldsymbol{n} + o(\delta)$, where $\boldsymbol{n}=(n,0) \in T_{\m}\mathbb{S}^2$, into \eqref{e:llg} gives, at order $\delta$, the linear equation 
\[
\partial_t \boldsymbol{n} =(\pm\mu-h) \boldsymbol{n} \times \3 
\pm \3\times  (\alpha \partial_t \boldsymbol{n} - \partial_x^2 \boldsymbol{n} + {\beta^\pm}  \boldsymbol{n}\times \3),
\]
which may be written in complex form as
\[
\partial_t n \pm {\beta^\pm} n = i \left( \alpha \partial_t n - \partial_x^2 n - (h \mp \mu) n \right).
\]
Its eigenvalue problem diagonalizes in Fourier space (for $x$) and yields the matrix eigenvalue problem
\[
\left | 
\begin{array}{ccc}
\pm{\beta^\pm} -\lambda & \Lambda \\
-\Lambda & \pm{\beta^\pm} -\lambda 
\end{array}
\right | = 0,
\]
where $\Lambda = \pm\mu - h \mp \alpha \lambda \mp \ell^2$ with $\ell$ the Fourier wave number.
The determinant reads 
\[
(\pm{\beta^\pm}- \lambda)^2 = - \Lambda^2 \quad \Leftrightarrow \quad
\pm{\beta^\pm} -\lambda = \sigma \rmi \Lambda, \; \sigma\in\{\pm 1\}.
\]
Considering real and imaginary parts this leads to
\begin{align*}
(1+\alpha^2)\Re(\lambda) &= \pm{\beta^\pm} - \alpha(\ell^2\pm h-\mu)= \alpha(\mu \mp(h-{\beta^\pm}/\alpha)-\ell^2) \\
\Im(\lambda) &= \sigma(\mp\Re(\lambda) +{\beta^\pm}/\alpha), 
\end{align*}
so that the maximal real part has $\ell=0$. At criticality, where $\Re(\lambda)=0$ the imaginary parts are $\pm{\beta^\pm}/\alpha$, which (if nonzero)  corresponds to a so-called Hopf-instability of the (purely essential) spectrum and we expect the emergence of oscillating solutions whose frequency at onset is ${\beta^\pm}/\alpha$ \cite{RaSch}. Since the critical mode has $\ell=0$, the onset of instability coincides with the aforementioned Hopf-bifurcation of the ODE associated with diffusionless \eqref{e:llg}.

The formulas for real- and imaginary parts immediately give the results mentioned in \S\ref{s:intro} and

\begin{Lemma}\label{l:conststab}
The constant state $\m=\3$ is (strictly) $L^2$-stable if and only if $\mu < \mu^+:=h-{\beta^+}/\alpha$ and $\m=-\3$ if and only if $\mu <\mu^-:=-(h-{\beta^-}/\alpha)$. Instabilities are of Hopf-type for the essential spectrum and have frequency $\beta/\alpha$.

For $\cc=0$ the anisotropy is subsubcritical precisely for $-\mu>|h-\beta/\alpha|$, and supercritical precisely for $\mu>|h-\beta/\alpha|$. 
\end{Lemma}

\section{Wavetrains}\label{s:wt}

To exploit the rotation symmetry about the $\3$-axis, we change to polar coordinates in the planar components $m=(m_1,m_2)$ of the magnetization
$\m = (m, m_3)$. With $m=r \exp(\rmi \varphi)$  equation \eqref{e:llg} changes to

\begin{equation}\label{e:llg-mod}
\left(\begin{array}{cc}
\alpha & -1\\
1 & \alpha
\end{array}\right)
\left(\begin{array}{c}
r^2 \partial_t \varphi\\
\partial_t m_3 
\end{array}\right)
=
\left(\begin{array}{c}
\partial_x(r^2\partial_x\varphi) \\
\partial_x^2 m_3 + |\partial_x \m|^2 m_3
\end{array}\right)
+
r^2
\left(\begin{array}{c}
\beta/(1+\cc m_3)\\
h-\mu m_3
\end{array}\right),
\end{equation}
where 
\[
|\partial_x \m|^2 = (\partial_x r)^2 + r^2 (\partial_x \varphi)^2 + (\partial_x m_3)^2 \quad \text{and}  \quad r^2 + m_3^2 =1.
\] 
This can be seen as follows. In view of \eqref{eq:LL_form}, with $\boldsymbol{h}_{\rm eff}$ replaced by the extended effective field 
$\boldsymbol{H}_{\rm eff}= \boldsymbol{h}_{\rm eff} - \m \times \boldsymbol{j}$ as in \eqref{eq:H_ext} and taking into account \eqref{eq:div}, 
\eqref{e:llg} reads
\[
\alpha \partial_t \m + \m\times  \partial_t \m = \partial_x^2 \m +(h- \mu m_3) \boldsymbol{\hat e}_3 + \left( |\partial_x \m|^2 + \mu m_3^2 - h m_3 \right) \m
- \frac{\beta}{1+\cc m_3} \m\times  \3.
\]
The third component of the above equation is the second component of \eqref{e:llg-mod}, whereas the first component of \eqref{e:llg-mod} is obtained upon inner multiplication by $\m^\perp = (m^\perp,0)=(ie^{i \varphi},0)$ and taking into account that  $\m \times \3= - \m^\perp$. 

The rotation symmetry has turned into the shift symmetry $\varphi \mapsto \varphi$ + const.
In full spherical coordinates $\m=\binom{e^{i \varphi} \sin \theta}{\cos \theta}$,  \eqref{e:llg-mod} further simplifies to 
\begin{equation}\label{e:llg-cyl}
\left(\!\begin{array}{cc}
\alpha & -1\\
1 & \alpha
\end{array}\!\right)
\left(\!\begin{array}{c}
\sin \theta \partial_t \varphi\\
-\partial_t \theta 
\end{array}\!\right)
=
\left(\!\begin{array}{c}
2\cos \theta \partial_x \theta \, \partial_x \varphi + \sin \theta\partial_x^2\varphi \\ 
-\partial_x^2 \theta + \sin\theta\cos\theta (\partial_x\varphi)^2
\end{array}\!\right)
+
\sin \theta
\left(\!\begin{array}{c}
\beta/(1+\cc \cos(\theta))\\
h-\mu \cos \theta 
\end{array}\!\right).
\end{equation}

\subsection{Existence of wavetrains}\label{s:ex-wt}

Wavetrains are solutions of the form $\m(x,t) = \m_*(kx-\omega t)$, where $k$ is referred to as the wavenumber and $\omega$ as the frequency. A natural type of wavetrains are relative equilibria with respect to the phase shift symmetry for which
$\varphi = k x - \omega t$ and $m_3, r$ are constant. See Figure~\ref{f:wavetrain} for an illustration.

\begin{Theorem}\label{t:wt-ex-gen} Wavetrains with frequency $\omega$ and wavenumber $k$ are in one-to-one correspondence to solutions of
 \[
\omk(\omega,k):=\cc\alpha\omega(\omega+h)-(\beta+\alpha \omega)(k^2-\mu) =0,
\]
under the constraint $|(\omega+h)/(k^2-\mu)|\leq 1$. In particular, for each $k$ there are at most two values of $\omega$ that yield a wavetrain, and for each $\omega\neq -\beta/\alpha$ there is at most one value of $k^2$ that gives a wavetrain, unless $\cc\alpha\omega(\omega+h)=0$ for $\omega=-\beta/\alpha$. Moreover, for $|\cc|<1$, 

\begin{itemize}
\item[(a)] $\omega\neq0$, $k^2\neq \mu$, $\sgn(\omega)=-\sgn(\beta)$ and $\omega\in [\min\{-\beta^\pm/\alpha\}, \max\{-\beta^\pm/\alpha\}]$.
\item[(b)] As $|k|\to\infty$ we have $\omega \to -\beta/\alpha$ and $m_3\to 0$.
\item [(c)] Bifurcations of $k\sim 0$ from $k=0$ for fixed $\omega$ are unfolded for increasing $\mu$.
\item[(d)] Bifurcations for $\m\neq \pm\3$ occur at $\omega=\omega_{\rm sn} := \frac{\beta\pm\sqrt{\beta(\beta-4\alpha h)}}{2\alpha}$ with $|m_3|<1$ if $|\frac{1}{\cc} - \frac{h}{\mu}|<2$ and are either:

... a hyperbolic point for $k\neq0$, if $h=\beta/\alpha$ and then $\omega=-h, k^2=-\cc h+\mu$,

... a hyperbolic point for $k=0$, if $\sgn((\cc h-\mu)(\cc h+\mu))=1$,

... an elliptic point for $k=0$, if $\sgn((\cc h-\mu)(\cc h+\mu))=-1$.
\end{itemize}
\end{Theorem}

\begin{figure}[t]
\begin{center}
\begin{tabular}{cc}
\includegraphics[width=0.3\textwidth]{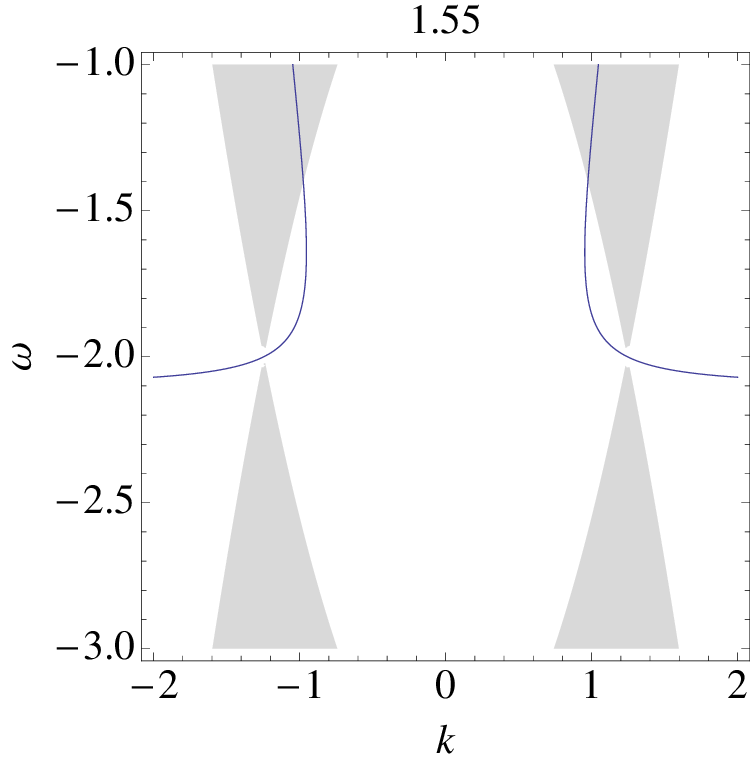}
& \includegraphics[width=0.3\textwidth]{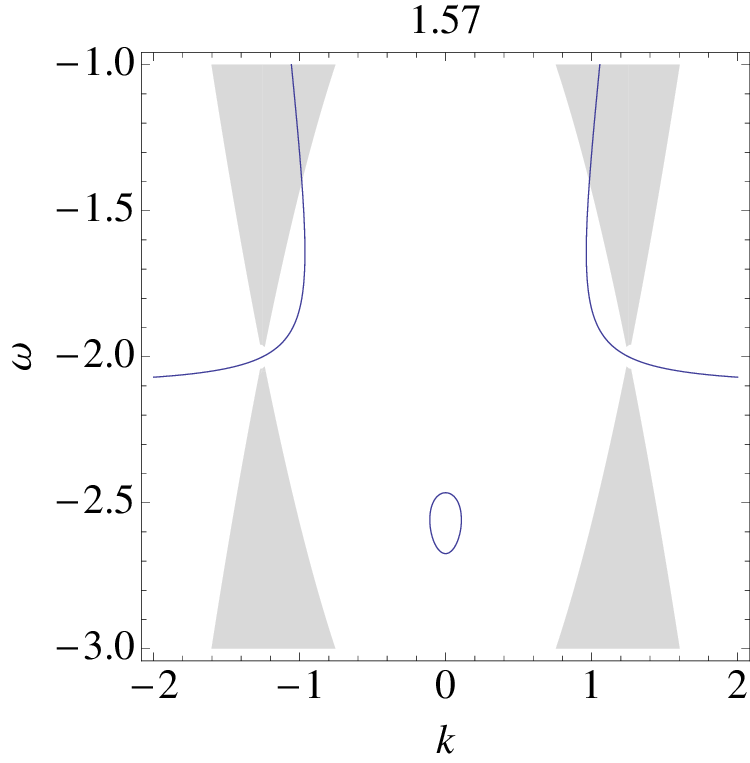}
\end{tabular}
\caption{Plots of wavetrain locations in the $(k,\omega)$-plane when parameters pass through an elliptic bifurcation point. Labels are the values of $\mu$. Other parameters are fixed at $\cc=0.5, h=2, \alpha=1, \beta=2.1$ so that $\omega^\pm=-\beta^\pm$, $\beta^+=1.4, \beta^-=4.2$ and $\mu^+=0.6, \mu^-=2.2$. Therefore, $\mu=1.55, 1.57$ are both subcritical with $\3$ stable and $-\3$ unstable. Shaded regions have $|m_3|>1$. The elliptic point lies at $\mu_{\rm sn} \approx 1.558, \omega_{\rm sn} = -2.558$.}
\label{f:wt-elliptic}
\end{center}
\end{figure}

\begin{figure}[t]
\begin{center}
\begin{tabular}{ccc}
\includegraphics[width=0.3\textwidth]{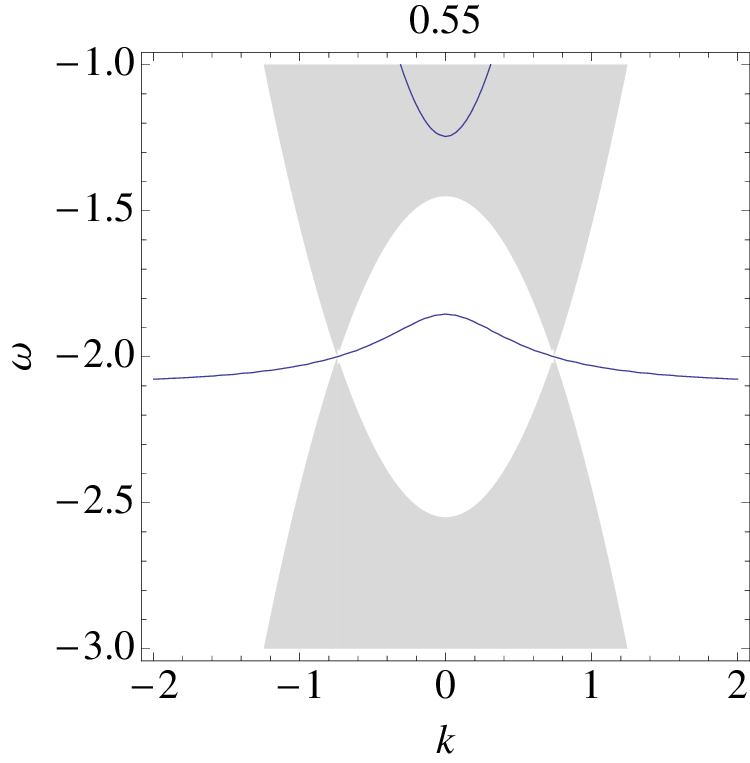}
& \includegraphics[width=0.3\textwidth]{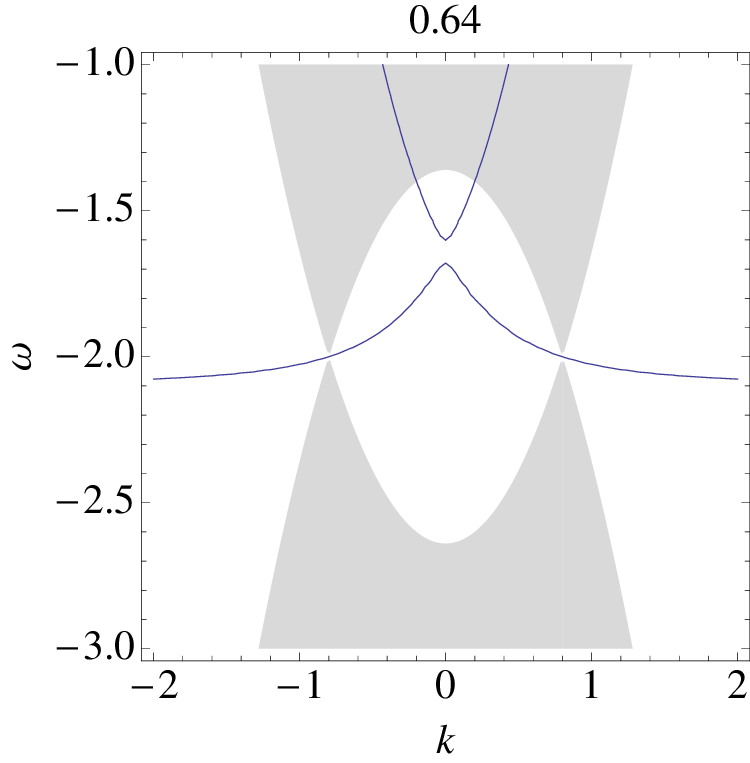}
& \includegraphics[width=0.3\textwidth]{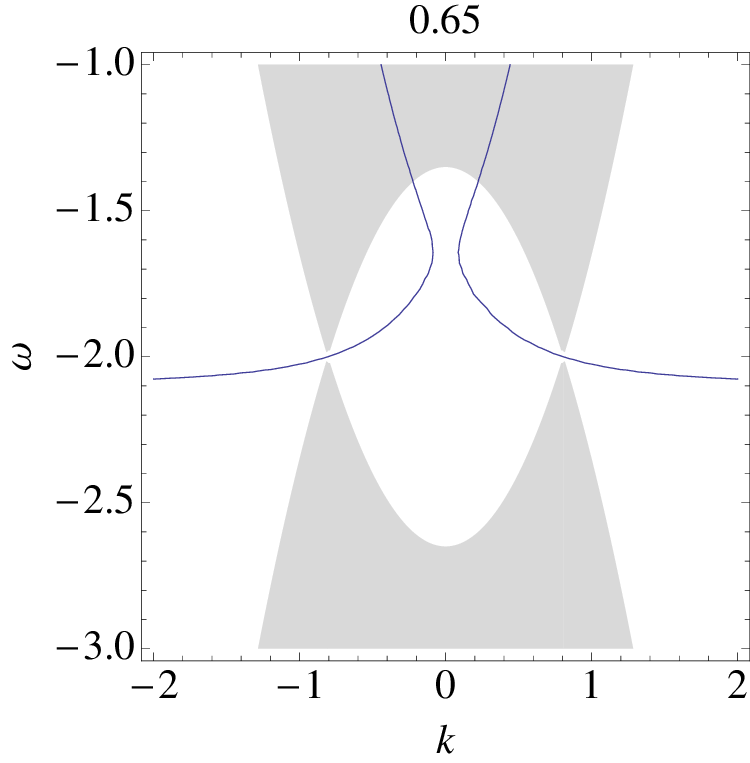}\end{tabular}
\caption{Analogue of Figure~\ref{f:wt-elliptic} with fixed parameters as there, when $\mu$ passes through a hyperbolic bifurcation point at $k=0$ with $\mu_{\rm sn} \approx 0.641, \omega_{\rm sn} = -1.641$. Note that between $\mu=0.55$ and $\mu=0.64$ the upper branch enters the region $|m_3|\leq 1$ at a bifurcation of in this case $\3$ at $\mu=\mu^+=0.6$. Therefore, $\mu=0.64, 0.65$ are supercritical with $\pm\3$ both unstable.}
\label{f:wt-hyper}
\end{center}
\end{figure}

\begin{figure}[t]
\begin{center}
\begin{tabular}{ccc}
\includegraphics[width=0.3\textwidth]{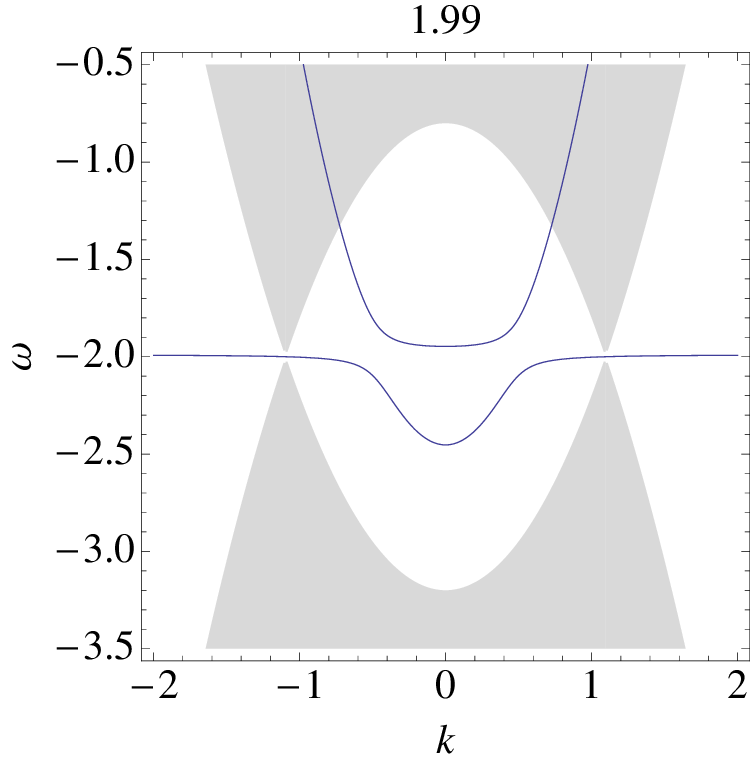}
& \includegraphics[width=0.3\textwidth]{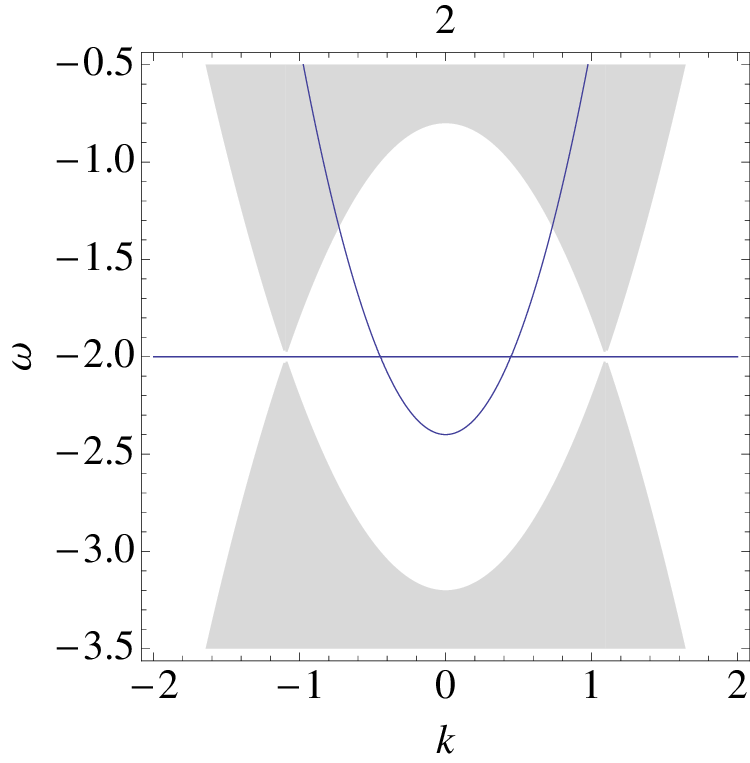}
& \includegraphics[width=0.3\textwidth]{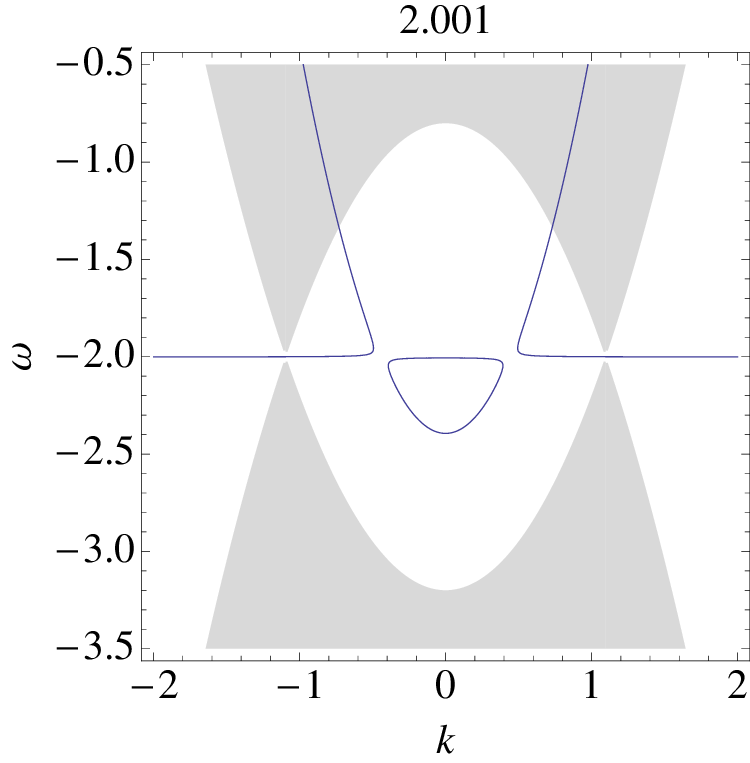}\\
\end{tabular}
\caption{Analogue of Figure~\ref{f:wt-elliptic} when parameters pass through hyperbolic bifurcation points with $k\neq0$ so that $h=\beta/\alpha$. Notably this cannot be unfolded by variation of $\mu$. Labels are the values of $\beta$ and other parameters are fixed at $\cc=0.5, h=2, \alpha=1, \mu=1.2$. This is subcritical with $\3$ unstable since $\mu^+=2(1-\beta/3), \mu^-=2(\beta-1)$ and $1.2\in(\mu^+,\mu^-)$. }
\label{f:wt-hyperB}
\end{center}
\end{figure}

In the following we discuss the existence problem and thereby prove each statement of the theorem. 

The terms elliptic and hyperbolic refers to the use in \cite{RaSch07} and will be explained below. Note that the latter two bifurcation types do not occur for $\cc=0$. That case is considered in detail in \S\ref{s:exc0}. We note that $\Gamma$ depends on $k$ only through $k^2-\mu$ so that $\mu\geq0$ is the same as $\mu=0$ up to change in wavenumber and solutions for fixed $\omega$ can increase $|k|$ only through increasing $\mu$.

Substituting the wavetrain ansatz into \eqref{e:llg-cyl} yields the algebraic equations
\[
\left(\begin{array}{cc}
\alpha & -1\\
1 & \alpha
\end{array}\right)
\left(\begin{array}{c}
-\sin(\theta) \omega\\
0
\end{array}\right)
=
\left(\begin{array}{c}
0\\
\sin(\theta)\cos(\theta)k^2
\end{array}\right)
+
\sin(\theta)
\left(\begin{array}{c}
\beta/(1+\cc \cos(\theta))\\
h-\mu \cos(\theta)
\end{array}\right).
\]
Thus either $\theta\equiv0 \mod \pi$ or (recall $|\cc|<1$)
\begin{equation}\label{e:wt-om}
-\alpha\omega = \frac{\beta}{1+\cc m_3}\;,\qquad -\omega  = (k^2-\mu)m_3 + h.
\end{equation}
In the first case we have $r=0$, which corresponds to the constant upward or downward magnetizations, $(r,m_3)=(0,\pm1)$ with unspecified $k$ and $\omega$. In the second case, we  notice aside that absence of dissipation ($\alpha=0$) requires absence of current ($\beta=0$) and there is a two-parameter set of wavetrains given by the second equation. The case we are interested in is $\alpha>0$ and then $\beta=0$ requires $\omega=0$, and this falls into the special case $\cc=0$.

In the generic case $\beta,\cc\neq0$, the first equation implies that $\omega\approx 0$ is not possible for $|m_3|\leq 1$ and we obtain
\begin{equation}\label{e:wt-m3}
m_3 = -\frac{1}{\cc}\left(\frac \beta{\alpha\omega}+1\right)\;, \qquad
m_3 = -\frac{\omega+h}{k^2-\mu}
\end{equation}
where for $k^2=\mu$ we have $\omega=-h$ and the first equation holds.

Eliminating $m_3$ and rearranging terms gives the existence condition in terms of $\omega$ and $k$ as zeros of $\omk(\omega,k)$ as in the theorem. For $\omega\neq-\beta/\alpha$ this gives $k^2$ as a quadratic function of $\omega$ inverse to the nonlinear dispersion relation $\omega(k)$. The exceptional $\omega=-\beta/\alpha$ occurs precisely when $h=\beta/\alpha$ and implies $m_3=0$, i.e.\ a solution on the equator (or $\cc=0$). The strict monotonicity of $\omk$ in $k$ away from $k=0$ also means that upon parameter change new solution branches can emerge only through local extrema of $\omk$ at $k=0$, i.e., an `elliptic' point. Specifically, this occurs if at a critical point $\partial_k^2\omk(\omega,0)=2(\beta+\alpha\omega)$ has the same sign as $\partial_\omega^2\omk(\omega,0)=2\alpha\cc$ and, e.g. $\omk(0,0)=\beta\mu$ varies. 

In particular, $\omk(\omega,k)=\partial_\omega\omk(\omega,k)=0$ occurs at 
\begin{align}\label{e:omcrit}
\alpha\omega^2 + \beta(2\omega+h)=0
\end{align}
and is a fold point of wavetrains (in the form of homogeneous oscillations) with fixed $k$. For $k=0$ these have at frequency and parameters (recall $\cc\neq0$)
\begin{align}\label{e:omsn}
\omega_{\rm sn} = -\frac{\cc h+\mu}{2\cc}\;, \qquad 4\beta\mu\cc =\alpha (\mu+\cc h)^2.
\end{align}
Note that $|m_3|<1$ for $k=0$ is by \eqref{e:wt-m3} equivalent to $|\omega_{\rm sn}+h|<|\mu|$ which yields $|\frac 1 \cc - \frac h \mu|<2$.

At such critical point we also have
\begin{align}\label{e:crit}
\partial_k^2\Gamma(\omega,0) &= -\frac{\alpha}{\cc}(\cc h-\mu)(\cc h+\mu),
\end{align}
so that the relative size of $\cc h$ and $\mu$ determines whether such a bifurcation point is elliptic or hyperbolic in the language of \cite{RaSch}.
In terms of critical parameters, substituting $\omega=\omega_{\rm sn}+\tomega$ and, for instance $\mu=\mu_{\rm sn}+\tmu$ to unfold with $\mu$ and other parameters fixed we obtain
\begin{align}\label{e:wt-hypell}
\tmu = k^2- \frac{\cc\alpha}{\beta+\alpha\omega_{\rm sn}}\tomega^2
\end{align}
which gives the options of hyperbola or ellipses for level sets. Hyperbolic points are saddle points of $\omk$ and at such points the connectivity of existing branches changes. For $k\neq0$ this occurs in particular, if $h=\beta/\alpha$ when the two branches of $\omk=0$ are the line $\omega=-h$ for any $k$ and $\omega = \frac{1}{\cc}(k^2-\mu)$.

We plot examples of these situations in Figures~\ref{f:wt-elliptic}, \ref{f:wt-hyper}, \ref{f:wt-hyperB}.

More globally, since $\omk$ is quadratic in $\omega$ there are at most two solutions for each $k$ and by strict monotonicity in $k^2$, away from $h=\beta/\alpha$, there is at most one solution for each $\omega$ or the whole line $\omega=-\beta/\alpha$. The only complication is the constraint $|m_3|\leq 1$ -- dispersion curves touch the boundary $m_3=1$ at bifurcations of $\pm\3$, which were studied in \S\ref{s:hopf}. The figures illustrate the essential scenarios. 


\subsection{Existence in the case $\cc=0$}\label{s:exc0}

\begin{figure}
\begin{center}
\begin{tabular}{ccc}
\includegraphics[width=0.3\textwidth]{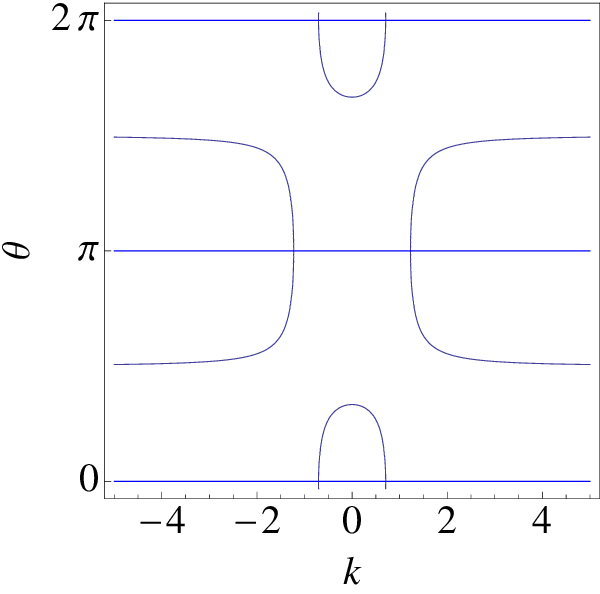}
& \includegraphics[width=0.3\textwidth]{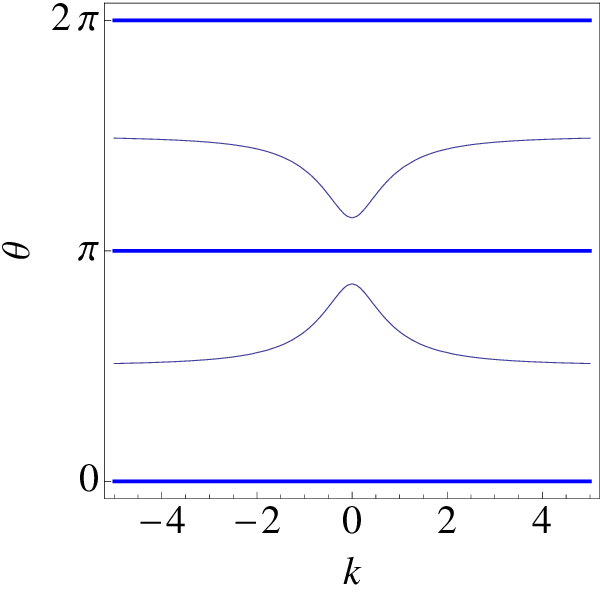}
& \includegraphics[width=0.3\textwidth]{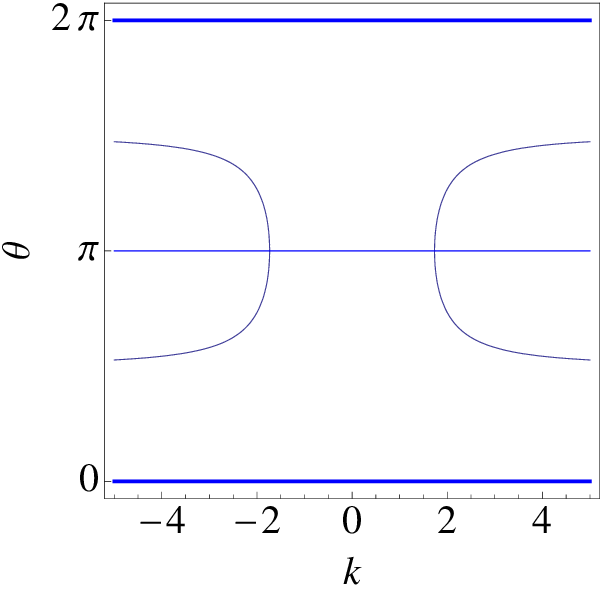}\\
(a) & (b) & (c)
\end{tabular}
\caption{Plots of equilibrium locations in the $(k,m_3)$-plane including the trivial equilibria $\pm\3$ plotted with thick line if stable (when not intersecting wavetrain parameters). Compare Figure~\ref{f:wt-ex}. See \eqref{e:psi}. (a) supercritical (easy plane) anisotropy ($\mu=1$, $h-\beta/\alpha=1/2$), (b) subcritical anisotropy ($\mu=1$, $h-\beta/\alpha=2$), where no homogeneous oscillations ($k=0$) exist, and (c) subsubcritical (easy axis) anisotropy $\mu=-1$, $h-\beta/\alpha=0.9$.}
\label{f:wavetrains}
\end{center}
\end{figure}

In this case the existence conditions can be conveniently written as 
\begin{align}
\omega &= -\frac{\beta}{\alpha} \label{e:freq}\\
\cos(\theta) &= \frac{h-\beta/\alpha}{\mu-k^2}, \quad (\mu\neq k^2). \label{e:psi}
\end{align}
As expected from the variational structure in rotating coordinates discussed in \S\ref{s:symc0}, all wavetrains oscillate with frequency given by the ratio of applied current and dissipation. In particular, in this case the natural representation of wavetrains is that $(\theta,k)$-plane rather than $(\omega,k)$ as above.

An involution symmetry involving parameters is 
\begin{equation}\label{e:parasym}
(h-\beta/\alpha, \theta) \to (\beta/\alpha-h, \theta + \pi),
\end{equation} 
so that the sign of $h-\beta/\alpha$ is irrelevant for the qualitative picture.

Solvability of \eqref{e:psi} requires that $|\mu-k^2| > |h-\beta/\alpha|$ (unless $r=0$), so that only for super- and subsubcritical anisotropy, 
\begin{equation}\label{e:muex}
|\mu|>|h-\beta/\alpha|,
\end{equation}
there exist wavetrains with wavenumber in an interval around $k=0$. In other words, non-trivial spatially nearly homogeneous oscillations require sufficiently small (in absolute value) difference between applied magnetic field and oscillation frequency (ratio of applied current and dissipation). The transition into this regime goes via the `Hopf' instability from \S\ref{s:hopf}. Combining \eqref{e:muex} with Lemma~\ref{l:conststab} and straightforward analysis of \eqref{e:psi} gives the following lemma. The three types of solution sets are plotted in Figure~\ref{f:wavetrains}. Clearly, the solution sets are symmetric with respect to the signs of $k$ and $\theta$, respectively.

\begin{figure}
\begin{center}
\scalebox{0.8}{\input{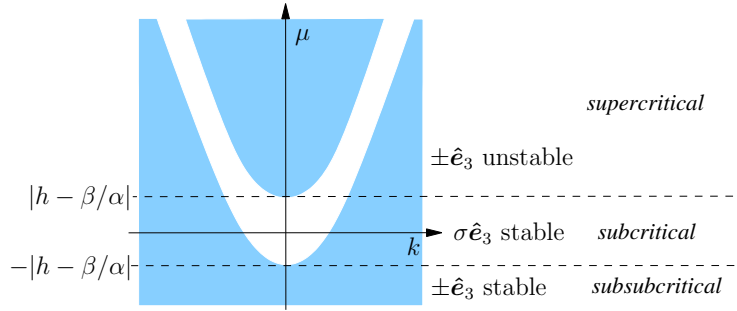}}\\
\caption{Sketch of existence region (shaded) in the $(k,\mu)$-plane, with boundary given by \eqref{e:muex} for $\sigma:=\mathrm{sgn}(h- \beta/\alpha)\neq 0$. The sign of $\sigma$ determines which of $\pm\3$ is stable in the subcritical range.}
\label{f:wt-ex}
\end{center}
\end{figure}

\begin{Theorem}~\label{t:wt-ex}
There are three types of wavetrain parameter sets solving \eqref{e:psi}: 
\begin{enumerate}
\item For supercritical anisotropy there is one connected component of wavetrain parameters including $k=0$, and two connected components with unbounded $|k|$, each with constant sign of $k$.
\item For subcritical anisotropy there are two connected components with unbounded $|k|$, each with constant sign of $k$.
\item For subsubcritical anisotropy there are two connected components, each a graph over the $k$-axis.
\end{enumerate}

The Hopf-type instabilities of $\pm\3$ noted in \S\ref{s:hopf} at the transition from sub- to supercritical anisotropy is a supercritical bifurcation in the sense that solutions emerge at the loss of stability of the basic solution, here $\pm\3$, while that from sub- to subsubcritical is subcritical in the sense that solutions emerge at the gain of stability.
\end{Theorem}

\begin{Proof}
Let us consider the existence region of wavetrains in wavenumber-parameter space. From \eqref{e:psi}, $r=0$ at $\theta\equiv 0 \mod \pi$ gives the boundary for nontrivial amplitude, 
\begin{align}\label{e:exbdry}
\mu= k^2 \pm(h-\beta/\alpha),
\end{align}
as a pair quadratic parabolas in $(k,\mu)$-space. The solution set in this projection is sketched in Figure~\ref{f:wt-ex}. Remark that this set is non-empty for any parameter set $\alpha, \beta, h\in \R$ of \eqref{e:llg}. However, as in the general case not all wavenumbers are possible due to the geometric constraint. Notably, the existence region consists of two disjoint sets, one contained in $\{\mu>0\}$ with convex boundary and one extending into $\{\mu<0\}$ with concave boundary. 
\end{Proof}

\subsection{Stability of wavetrains}\label{sec:spectral}

In this section we discuss spectral stability of wavetrains and in summary we obtain the following result.

\begin{Theorem}\label{t:wt-stab}~
\begin{itemize}
\item[(a)] Wavetrains bifurcating from $\pm\3$ at $k=0$ are stable if $\mu^\pm >\max\{0,\cc\beta^\pm/\alpha\}$, which implies supercritical bifurcation, i.e., for increasing $\mu$. They are unstable if $\mu^\pm <\max\{0,\cc\beta^\pm/\alpha\}$, which means subcritical bifurcation, i.e., decreasing $\mu$.
\item[(b)] For $\mu>0$ there is $k_*>0$ such that precisely the wavetrains with wavenumber $|k|<k_*$ and $\mu>k^2+\alpha\omega^2\cc/\beta$ are sideband stable. These are fully spectrally stable if $\beta\cc\geq 0$ or $\alpha\geq0$ is sufficiently small. A wavetrain and $\3$ or $-\3$ can be simultaneously stable only for $\cc\neq 0$.
\item[(c)] For each $k$ at most one wavetrain can be stable. Wavetrains with $r\sim 0$ and $k\neq0$ are unstable. All wavetrains with $k^2>\mu$ are unstable. 
\item[(d)] Near the hyperbolic or elliptic bifurcation points at $k=0$ the sideband stable wavetrains lie in a sector that is to leading order bounded by $|\omega-\omega_{\rm sn}|=S|k|$ whose opening angle less than $\pi$ and which includes $\omega$ with a selected sign of $\omega-\omega_{\rm sn}$.
\end{itemize}
\end{Theorem}

\begin{Remark}
Item (a) should be compared with the Hopf instabilities discussed in \S\ref{s:hopf}. The case $\cc=0$ is simplest: the constraint simply means that wavetrains bifurcating at the subsub- to subcritical transition are unstable. For general $\cc$ the condition also accounts for interaction with folds.

Concerning item (b), notably \emph{in the easy axis case} $\mu<0$ all wavetrains are unstable. The condition $\cc\beta\geq0$ is not sharp. However, we do not know whether a `Hopf' instability can occur for otherwise stable wavetrains if $\cc\beta<0$ is sufficiently large negative. 
Item (d) is an analog to the results in \cite{RaSch07}. 
\end{Remark}

\begin{figure}
\begin{center}
\begin{tabular}{cccc}
\includegraphics[width=0.22\textwidth]{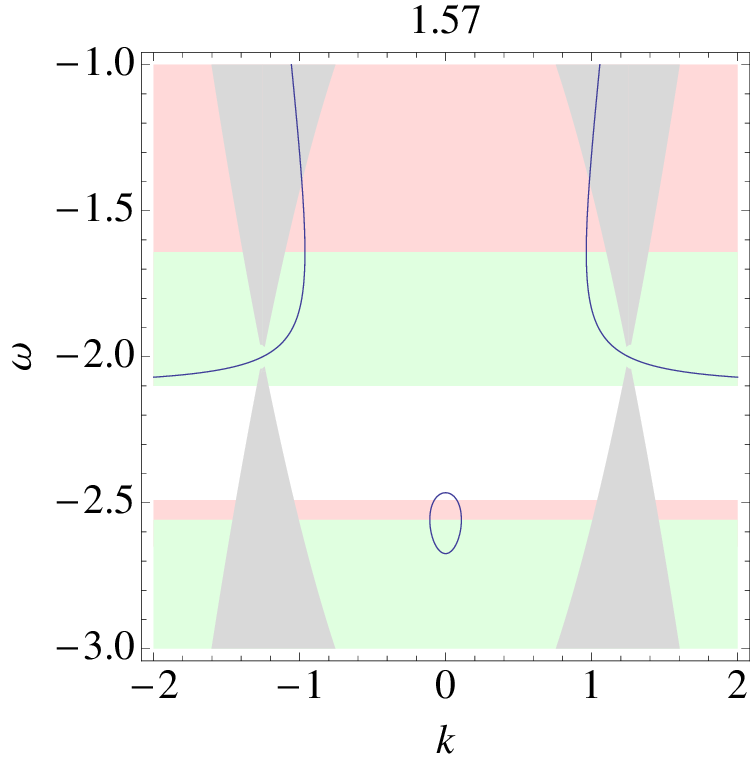}
& \includegraphics[width=0.22\textwidth]{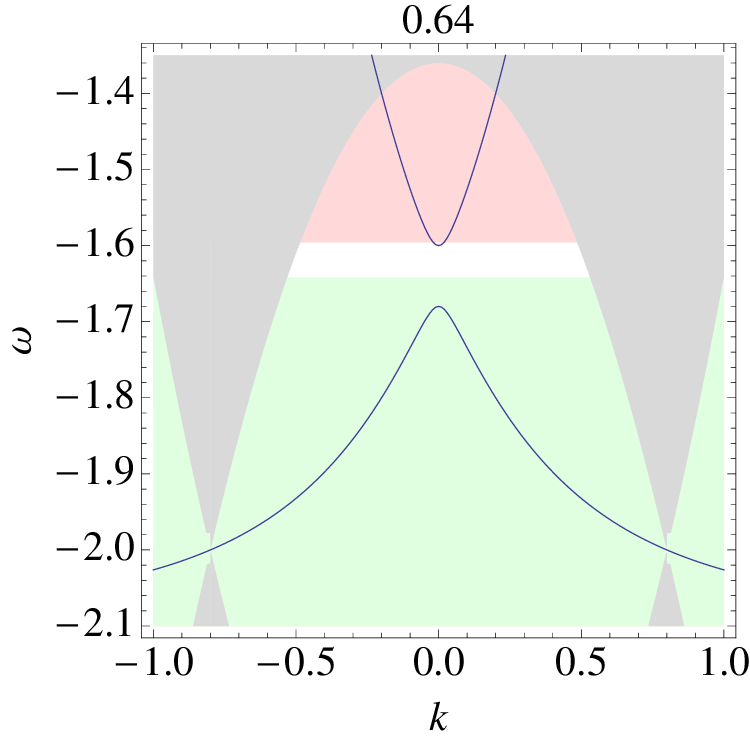}
& \includegraphics[width=0.22\textwidth]{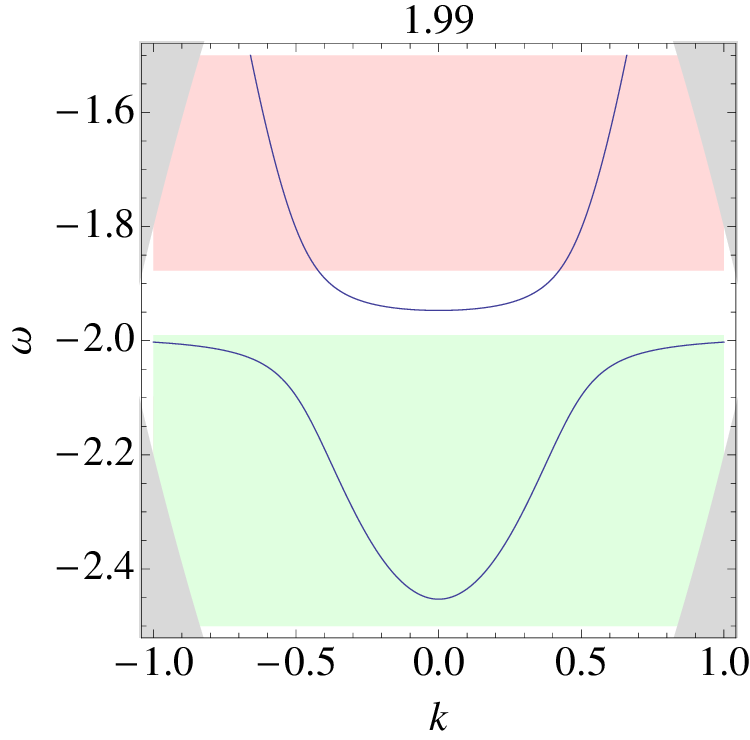}
& \includegraphics[width=0.22\textwidth]{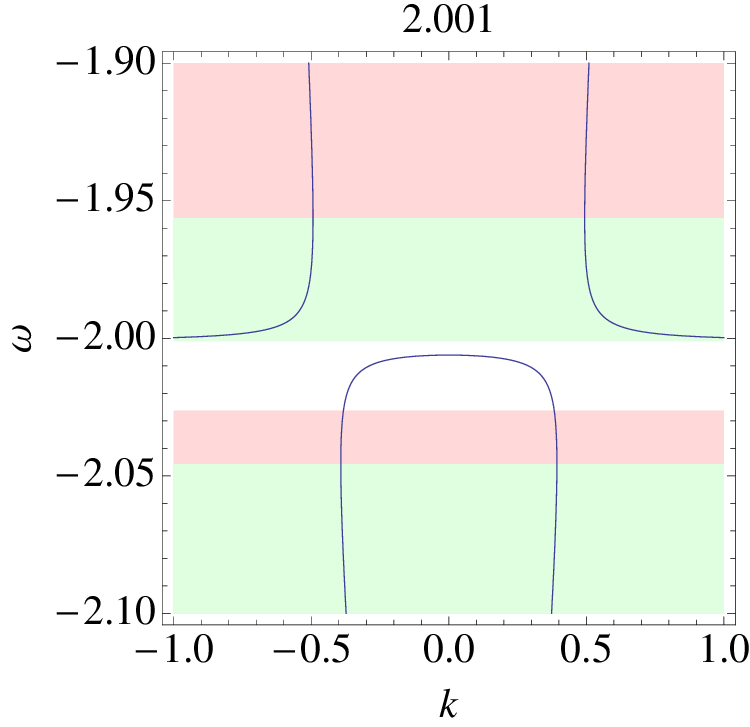}
\end{tabular}
\caption{Dispersion curves and stability regions. Wavetrains on the dispersion curves in the green shaded region are unstable due to unstable eigenvalue of $A(0,0)$ using \eqref{e:foldind}, and in the red shaded region due to unstable sideband only, using $\tilde\lambda_0''(0)$ with $k^2$ from $\omk=0$. Parameters are as in (a): Figure~\ref{f:wt-elliptic}; (b): Figure~\ref{f:wt-hyper}; (c),(d): Figure~\ref{f:wt-hyperB}. Only the wavetrains near $k=0$ are stable (except for the lower branch in (b)).}
\label{f:wt-stab}
\end{center}
\end{figure}

\bigskip
In order to study spectral stability of wavetrains we consider the comoving frame $y=x-\cph t$ with $\cph=\omega/k$ the wavespeed so that the wavetrain is an equilibrium of \eqref{e:llg-mod}. For convenience, time is rescaled to $t = (1+\alpha^2)\tt$. The explicit formulation of \eqref{e:llg-mod} then reads
\begin{equation}\label{e:llg-mod-ex}
\partial_\tt\left(\!\!\begin{array}{c}
\!  \varphi \!\\
\! m_3 \!
\end{array}\!\!\right)
=
\left(\!\!\begin{array}{c}
\alpha(\partial_y(r^2\partial_y\varphi)/r^2 +\tbeta(m_3)) +  
(\partial_y^2 m_3 + |\partial_y \m|^2 m_3)/r^2 + h - \mu m_3 + \cph \partial_y\varphi \\
\alpha(\partial_y^2 m_3 + |\partial_y \m|^2 m_3+r^2(h-\mu m_3)) - \partial_y(r^2\partial_y\varphi) - r^2\tbeta(m_3) + \cph\partial_y m_3
\end{array}\!\!\right),
\end{equation}
where $r^2 = 1-m_3^2$ and $\tbeta(m_3):=\beta/(1+\cc m_3)$.

Let $\calF=(\calF_1,\calF_2)^\rmt$ denote the right hand side of \eqref{e:llg-mod-ex}. Wavetrains have constant $r$ and $m_3$ so that quadratic terms in their derivatives can be discarded for the linearization $\calL$ of $\calF$ in a wavetrain, and from $|\partial_y \m|^2$ only $r^2(\partial_y \varphi)^2$ is relevant. The components of $\calL$ are
\begin{align*}
\partial_\varphi \calF_1 &= \alpha \partial_y^2 + \cph \partial_y + 2k m_3 \partial_y\\
\partial_{m_3} \calF_1 &= r^{-2}\partial_y^2 + k^2 - \mu - 2\alpha k m_3r^{-2} \partial_y + \alpha \tbeta'(m_3)\\
\partial_{\varphi} \calF_2 &= 2\alpha k m_3 r^2 \partial_y - r^2 \partial_y^2\\
\partial_{m_3} \calF_2 &= \alpha \partial_y^2 + \alpha r^2(k^2-\mu) + \cph \partial_y - 2\alpha m_3(m_3 k^2 + h - \mu m_3) + 2m_3(k\partial_y + \tbeta(m_3))\\ &\quad- r^2\tbeta'(m_3)
 ~=~  \alpha \partial_y^2 + \alpha r^2(k^2-\mu) + \cph \partial_y + 2m_3 k\partial_y - r^2\tbeta'(m_3),
\end{align*}
where the last equation is due to \eqref{e:wt-m3} and $\tbeta'(m_3)= -\cc\beta/(1+\cc m_3)^2$. Since all coefficients are constant, the eigenvalue problem 
\[
\calL u = \lambda u
\]
is solved by the characteristic equation arising from an exponential ansatz $u=\exp(\nu y)u_0$, which yields the matrix 
\begin{align*}
A(\nu, \cph) := 
\left(\begin{array}{cc}
\alpha \nu^2 + (\cph+2 k m_3)\nu & -r^{-2} \nu(-\nu + 2\alpha k m_3) + k^2-\mu+\alpha\tbeta'(m_3)\\
r^2\nu(-\nu + 2\alpha k m_3) & \alpha \nu^2 + (\cph + 2k m_3)\nu  + \alpha r^2(k^2-\mu)-r^2\tbeta'(m_3)
\end{array}\right).\\
\end{align*}
The characteristic equation then reads
\begin{align}
d_{\cph}(\lambda,\nu) & :=  |A(\nu,\cph) - \lambda| = |A(\nu,0) -(\lambda -  \nu \cph)| 
= d_0(\lambda-\cph \nu,\nu) \nonumber \\
d_0(\lambda,\nu) & = \lambda^2 - t(\nu) \lambda + d(\nu) \label{e:disp}
\end{align}
with trace and determinant of $A(\nu,0)$ 
\begin{align*}
t(\nu) &:=\tr A(\nu,0)  = 2\nu(\alpha \nu +2 k m_3)  + \alpha r^2(k^2-\mu)-r^2\tbeta'(m_3)\\
d(\nu) &:= \det A(\nu,0) = (1+\alpha ^2)\nu\left(\nu(\nu^2+ 4k^2 m_3^2 + r^2(k^2-\mu)) - 2r^2\tbeta'(m_3) km_3\right)\\
& = (1+\alpha ^2)\nu\left(\nu\left(\nu^2+ (3k^2 + \mu)\frac{(\tbeta(m_3)/\alpha - h)^2}{(k^2-\mu)^2} + k^2 - \mu\right) -2r^2 \tbeta'(m_3) km_3\right).
\end{align*}
In the last equation \eqref{e:psi} was used.

The characteristic equation is also referred to as the \emph{complex  (linear) dispersion relation}.
The spectrum of $\calL$, for instance in $L^2(\R)$, consists of solutions for $\nu=\rmi \ell$ and is purely essential spectrum (in the sense that $\lambda-\calL$ is not a Fredholm operator with index zero). Indeed, setting $\nu=\rmi\ell$ corresponds to Fourier transforming in $y$ with Fouriermode $\ell$. Note that the solution $d(0,0)=0$ stems from spatial translation symmetry in $y$. For the same reason the real part of solutions $\lambda$ of \eqref{e:disp} for any given $\nu\in\rmi\R$ does not depend on $\cph$, which means that spectral stability is independent of $\cph$ and is therefore completely determined by $d_0(\lambda, \rmi \ell)=0$.

As a first observation concerning stability we note that at $r=0, |m_3|=1$ the solutions $\lambda_\pm$ to $d_0(\lambda,\rmi\ell)=0$ have $\Re(\lambda_\pm) = -\ell(\ell\pm2k)$ whose maximum is $k^2$. Hence, wavetrains at (and thus near) $r=0$ are unstable for $k\neq 0$.

\bigskip
\paragraph{Fold stability.} First note the eigenvalues $A(0,0)$ are $0$ and 
\[
\tau(k^2):=t(0)=r^2(\alpha (k^2-\mu)-\tbeta'(m_3))
\]
so that a `fold instability' occurs precisely for $r=0$ (compare \S\ref{s:hopf}) or $\alpha(k^2-\mu)=\tbeta'(m_3)$. We readily compute that the latter corresponds to the critical points in \eqref{e:omcrit}. Note that 
$\tbeta'(m_3) = -\alpha^2\omega^2\cc/\beta$. In particular, for the loops of wavetrains emerging from an elliptic point, the upper and lower $\omega$-values have opposite signs of $\tau(0)$; specifically the lower is stable if $\beta\cc>0$, cf.\ Figure~\ref{f:wt-elliptic}. 

At the bifurcation points of $\pm\3$ (see \S\ref{s:hopf}) where $r=0$, we have $\mu=\mu_\pm$ and $\omega=-\beta^\pm/\alpha$ which yields Theorem~\ref{t:wt-stab}(a) except for the super-/subcriticality. To see this recall that $\pm\3$ destabilize always through increasing $\mu$. Whether wavetrains with $k=0$ emerge from $\pm\3$ depends on the sign of $\partial_\mu m_3$. From \eqref{e:wt-om} we find $\mu = \frac{h}{m_3} - \frac{\beta}{\alpha m_3(1+\cc m_3)}$ and compute 
\[
\partial_{m_3}\mu(\pm1) = \pm\left.\frac{\tau(0)}{\alpha r^2}\right|_{m_3=\pm 1}.
\]
Hence, fold-stability implies $\partial_\mu m_3(\mu^+)<0$ or $\partial_\mu m_3(\mu^-)>0$, respectively, so that increasing $\mu$ yields $|m_3|<1$ and thus emergence of solutions. 

More generally, sign changes of $\tau(k^2)$ correspond to fold points and a curve of spectrum crosses the origin. We plot some fold stability boundaries in Figure~\ref{f:wt-stab}. Using $\tau(k^2)$ we have that wavetrains with $\mu<k^2+\cc\frac{\alpha}{\beta}\omega^2$ are unstable. In particular, for $\cc \beta\geq 0$ all wavetrains with $k^2>\mu$ are unstable. Using the existence condition we may also write this condition independent of $k$ as
\begin{align}\label{e:foldind}
\frac{\omega(\beta(2\omega+h)+\alpha\omega^2)}{\beta(\beta+\alpha\omega)}<0,
\end{align}
which explains the changes in the fold stability indicator at $\omega=-\beta/\alpha$ in the figure.

Coming back to $k=0$, at the bifurcation points of $\pm\3$ (see \S\ref{s:hopf}) where $r=0$, we have $\mu=\mu_\pm$. Fold-stability is then $\mu>\alpha\cc \omega^2/\beta$ which holds if 
\[
\pm\left(h-\frac{\beta^\pm}{\alpha}\right) > \frac{\beta\cc}{\alpha(1\pm\cc)^2} 
\Leftrightarrow \pm h > \frac{\beta}{\alpha}\frac{2\cc\pm 1}{(1\pm\cc)^2},
\]
as noted in Theorem~\ref{t:wt-stab} item (a).

\medskip
We next check the other possible marginal stability configurations case by case.

\paragraph{Sideband instability.} A sideband instability occurs when the curvature of the curve of essential spectrum attached to the origin changes sign so that the essential spectrum extends into positive real parts. Let $\tilde\lambda_0(\ell)$ denote the curve of spectrum of $A(\rmi\ell,0)$ attached to the origin, that is $\tilde\lambda_0(0)=0$, and let  $'$ denote the differentiation with respect to $\ell$. Derivatives of $\tilde\lambda_0$ can be computed by implicit differentiation of $d_0(\lambda,\rmi \ell)=\lambda^2-t(\rmi \ell)+d(\rmi \ell)=0$.

This gives $\tilde\lambda_0'(0) = \rmi\frac{d'(0)}{t(0)}= -\rmi\tbeta'(m_3)\frac{2(1+\alpha^2)k m_3}{\alpha(k^2-\mu)-\tbeta'(m_3)}$ therefore
\[
\tilde\lambda_0''(0) = \left.-\frac{2d'(d'-t' t)+d'' t^2}{t^3}\right|_{\ell=0}=:\Lambda(k^2)\frac{2(1+\alpha^2)r^4}{\tau^3(k^2)},
\]
where the $\omega$-dependence is suppressed. Some calculations yield $\Lambda$ as a cubic polynomial in $K=k^2$ given by $\Lambda=a_3 K^3 + a_2 K^2 +  a_1 K + a_0$ with 
\begin{align*}
a_3&=-\alpha^2(4-3r^2),\\
a_2&=\alpha(2\tbeta'r^2+\alpha\mu(8-5r^2)),\\
a_1&=-\alpha^2\mu^2(4-r^2)-4(1-r^2)\alpha^2\tbeta'(m_3)^2 - (\tbeta'(m_3)^2+4\alpha\tbeta'(m_3)\mu)r^2,\\
a_0&=\mu (\tbeta'(m_3)+\alpha\mu)^2 r^2. 
\end{align*}
Notably, $a_0=\mu\tau(0)^2$ and also $a_3<0$ since $r\in(0,1)$. A wavetrain with wavenumber $k$ is therefore (strictly) \emph{sideband stable} precisely when $\Lambda(k^2)\tau(k^2)<0$ and \emph{sideband unstable} for $\Lambda(k^2)\tau(k^2)>0$. As mentioned above, $\tau(k^2)>0$ for large enough $k^2$ so that all such (already unstable) wavetrains are also sideband unstable as $a_3<0$ holds always. 

Since $\tau(k^2)>0$ is the unstable fold case, we next assume $\tau(k^2)<0$ so that  sideband stability is precisely $\Lambda(k^2)>0$. 

\medskip
\textit{Homogeneous oscillations.} For $k=0$ we obtain
$\tilde\lambda_0''(0)=\Lambda(0)\frac{2(1+\alpha^2)r^4}{\tau^3(0)}=-2(1+\alpha^2)\frac{\mu}{\tau(0)},
$
so that $\mu>0$ is required for sideband stability (given fold stability $\tau(0)<0$).



Let us study this situation near the Hopf instability of $\pm\3$, where $\mu\approx\mu_\pm = \pm(h-\beta^\pm/\alpha)$ with $\mu>\mu_\pm$ and $\omega = \beta^\pm/\alpha$. Hence, sideband stability of the emerging wavetrains requires  $h>\beta^+/\alpha$ or $h<\beta^-/\alpha$, that is, $\mu^\pm>0$. Compare Theorem~\ref{t:wt-stab}(a).

\medskip
\textit{Near homogeneous.}
To leading order $\Lambda(k^2)=0$ is $a_1k^2+a_0=\calO(k^4)$ and via \eqref{e:wt-hypell} we have $a_1=a_1(\tomega,k^2), a_0=a_0(\tomega,k^2)$, where at bifurcation $a_0(0)=\partial_\omega a_0(0)=0$. Upon expanding we therefore find sideband instabilities to leading order at
\[
\tomega^2 = \frac{a_1(0)-\partial_{k^2} a_0(0)}{\partial_{\tomega^2}a_0(0)}k^2 + \calO(|k|^3+|\tomega|^3),
\]
where $\partial_{\tomega^2}a_0(0) = \mu_{\rm sn} r_{\rm sn}^2\left(\frac{2}{\beta}\alpha^2 \omega_{\rm sn}\cc\right)^2$ and $\partial_{k^2} a_0(0) = \alpha \mu_{\rm sn}r_{\rm sn}^2$ with $r_{\rm sn}= 1-\left(\frac{\omega_{\rm sn}+h}{\mu}\right)^2$ the $r$-coordinate of the wavetrain at the bifurcation point. Some algebra yields $a_1(0) = -4\alpha(1+\alpha^2)(h+\omega_{\rm sn})^2$.

In accordance with the results in \cite{RaSch}, this means that wavetrains in a sector in the $(\omega,k)$-plane near the fold point are sideband stable, while wavetrains outside this sector are sideband unstable. We expect that the opening angle of this sector can be changed while keeping the dispersion curves essentially fixed. Here we do not pursue this further, but note that since sideband instabilities do occur and the stable region cannot include the fold points, the prefactor of $k^2$ is positive and $\tomega$ has a selected sign. The sideband boundaries for some examples are plotted in Figure~\ref{f:wt-stab}.

\medskip
\textit{General wavetrains.}
Concerning the sign of $\Lambda$ in general,  $\Lambda(0)=a_0$ has the sign of $\mu$ and thus all wavetrains for $\mu<0$ are sideband unstable.  For $\mu>0$ we have $\Lambda(\mu)=-4\alpha^2(\tbeta'(m_3))^2m_3^2\mu<0$, so a sign change occurs at some $K_{\rm sb}\in(0,\mu)$, which implies sideband instability for $k^2\geq \mu$. Moreover, $\mu>0$ and $\cc\beta\leq 0$ imply $a_1<0$ and $a_2>0$. Since $a_3<0$ this means both roots of $\Lambda'$, $(2a_2 \pm \sqrt{4(a_2^2-3a_3a_1)})/(6a_3)$ lie at negative $K$ and so $\Lambda$ is monotone decreasing for all $K>0$. Therefore $K_{\rm sb}$ is the unique sideband instability in this case.

On the other hand, for $a_1>0$ the roots of $\Lambda'$ have opposite signs so that due to the sign change in the interval $(0,\mu)$ the positive one must be a local maximum so that also in this case $K_{\rm sb}$ is the unique sideband instability. Moreover, in all cases wavetrains with $k^2>\mu$ are unstable since $\Lambda<0$ in this range.


\paragraph{Hopf instability.} A Hopf instability occurs when the essential spectrum touches the imaginary axis at nonzero values. In particular, there is $\gamma\neq0$ so that $d_0(\rmi\gamma,\rmi\ell)=0$. At $k=0$ we have $\Im(d_0(\rmi\gamma,\rmi\ell))=\gamma(2\alpha\ell^2-\tau(0))$ so that in the fold stable case $\tau(0)<0$ there is real $\ell$ for $\gamma\neq 0$. Recall $\tau(0)=-(\tbeta'(m_3) + \alpha\mu)r^2$. Therefore, there is no (relevant) Hopf instability near $k=0$.

More generally, solving $\Im(d_0(\rmi\gamma,\rmi\ell))$ for $\gamma$ and substituting the result into $\Re(d_0(\rmi\gamma,\rmi\ell))$ we obtain up to a factor $(1 + \alpha^2) \ell^2$ 
\[
\ell^2+(\mu-k^2)r^2 + 4\alpha^2m_3^2k^2 G(\ell^2), \; G(L)=-\frac{4L^2-4L(k^2-\mu)r^2+(\tbeta'(m_3)^2+(k^2-\mu)^2)r^4}{(2\alpha L -\tau(k^2))^2}
\]
where all terms except possibly $G(L)$ are positive in the interesting range $\mu>k^2$. Note that for sufficiently small $\alpha>0$ there is no root besides $\ell=0$ and thus no Hopf and in fact no sideband instability. But there seems to be no satisfying explicit bound. (While the same seems to occur for $m_3\sim0$, such wavetrains have $k^2>\mu$ and are thus unstable.)

However, $G$ is nondecreasing for $\tbeta'(m_3)\leq0$, i.e., $\cc\beta\geq 0$, since then
\[
G'(L) = -4\tbeta'(m_3) r^2\frac{2L+(\mu-k^2- \alpha \tbeta'(m_3))r^4}{(2\alpha L -\tau(k^2))^3}\geq 0,
\]
in the fold stable case $\tau(k^2)\leq 0$. Thus, besides $\ell=\gamma=0$, there is at most one solution $d_0(\rmi\gamma,\rmi\ell)=0$ for $\cc\beta\geq 0$, which rules out a Hopf instability as this requires two such solutions. We do not know whether or not Hopf instabilities can occur for general $\alpha$ and $\cc\beta< 0$. 

\paragraph{Turing instability.} A Turing instability occurs when the spectrum touches the origin for nonzero $\ell$, that is, there is $\ell\neq0$ so that $d_0(0,\rmi\ell)=0$, which means $\det A(\rmi\ell,0)=0$. Since $\Im(\det A(\rmi\ell,0))=-2(1+\alpha^2)\tbeta'(m_3)\ell k m_3r^2$ and our previous considerations already cover zeros of this, such instabilities do not occur.

\begin{figure}
\begin{center}
\begin{tabular}{cc}
\begin{minipage}{0.4\textwidth}
\hspace{5mm}\includegraphics[width=0.9\textwidth]{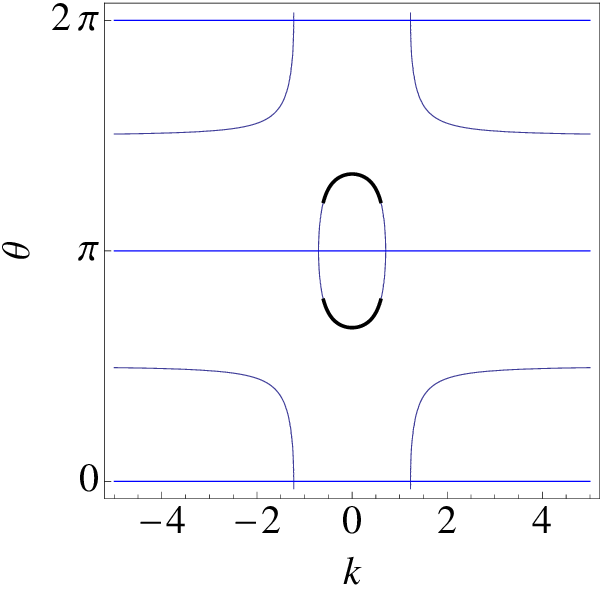}
\end{minipage} &
\hspace{5mm}\begin{minipage}{0.5\textwidth}
\scalebox{0.9}{\input{xfig/wt-stabregion.pstex_t}}\\ \\
\end{minipage}\\
(a) & (b)\hspace{10mm}
\end{tabular}
\caption{Analogues of Figures~\ref{f:wavetrains}(a) and~\ref{f:wt-ex} with stable range of wavetrains in (a) bold line, and in (b) the dark shaded region.}
\label{f:wavetrainstab}
\end{center}
\end{figure}

\bigskip
This exhausts the list of possible marginal stability. 

\subsection{Stability of wavetrains for $\cc=0$}

For $\cc=0$ the wavetrain frequency $\omega$ is independent of the wavenumber $k$ in \eqref{e:freq} so that the phase velocity $\omega/k$ of all wavetrains is $-\beta/(\alpha k)$ and the group velocity $\rmd \omega/\rmd k$, which describes the motion of perturbations by localized wave packets, vanishes for all wavetrains. 
Due to Theorem~\ref{t:wt-stab} destabilizations of stable wavetrains can only occur through a unique sideband instability $|k|=k_*>0$. This value can be explicitly determined since $d'(0)=0$ and $d''(0)=-2(1+\alpha^2)((3K+\mu)r^2-4K)$ so that $\tilde\lambda_0''(0)=0$ at 
\[
k_*^2=\mu \frac{r^2}{4-3r^2}.
\]
Taking into account that $k^2= \mu$ for a wavetrain can occur only if $h=\beta/\alpha$ we thus have

\begin{Theorem}\label{t:wt}
Consider $\cc=0$. All wavetrains whose wavenumber $k$ satisfies $|k|>k_*$ are unstable. For $\mu>0$ wavetrains with wavenumber $|k|<k_*$ are spectrally stable, while those with $|k|>k_*$ are unstable. In case $h\neq \beta/\alpha$ a sideband instability occurs at $k=\pm k_*$. There is no secondary instability for $k^2<\mu$. Nontrivial spectrally stable wavetrains exist only for supercritical anisotropy, $|h-\beta/\alpha|<\mu$.
\end{Theorem}

The overall picture for wavetrains of \eqref{e:llg} with $\cc=0$ is thus a combination of the scenarios from a supercritical and a subcritical real Ginzburg-Landau equation $\partial_t A = \partial_x^2 A + \tilde\mu A \mp A|A|^2, \; A(x,t) \in \C$, which describes the dynamics near pattern forming Turing instabilities and possesses the gauge-symmetry $A\to \rme^{\rmi\varphi}A$.  The interested reader is referred to the review \cite{AransonKramer:02} and the references therein.

\section{Coherent structures} \label{s:coh}

The coexistence of wavetrains and constant magnetizations raises the question how these interact. In this section we study solutions that have spatially varying local wave number. In particular, we consider solutions that spatially connect wavetrains or $\pm\3$ in a coherent way. In order to locate such solutions induced by symmetry we make the ansatz
\begin{equation}\label{e:ansatz}
\begin{array}{rl}
\xi &= x - st\\
\varphi &= \phi(\xi) + \Omega t\\
\theta &= \theta(\xi),
\end{array}
\end{equation}
with constant $s, \Omega$. Solutions of \eqref{e:llg-mod} of this form are generalized travelling waves to \eqref{e:llg} with speed $s$ that have a superimposed oscillation about $\3$ with frequency $\Omega$. This ansatz is completely analogous to that used in the aforementioned studies of the real and complex Ginzburg-Landau equations \cite{AransonKramer:02}. 

Let $\ttbeta(\theta):=\tbeta(\cos(\theta)) = \beta/(1+\cc \cos(\theta))$ denote the $(\beta,\cc)$-dependent term of  \eqref{e:llg-cyl}. Substituting ansatz \eqref{e:ansatz} into \eqref{e:llg-cyl} with $' = \rmd/\rmd \xi$ and $q=\phi'$ gives, after division by $\sin(\theta)$, the ODEs
\begin{equation}\label{e:llg-coh}
\left(\!\!\begin{array}{cc}
\alpha & -1\\
1 & \alpha
\end{array}\!\!\right)
\left(\!\!\begin{array}{c}
\sin(\theta) (\Omega-s q)\\
s \theta' 
\end{array}\!\!\right)
=
\left(\!\!\begin{array}{c}
2\cos(\theta)\theta' q + \sin(\theta)q' \\
-\theta'' + \sin(\theta)\cos(\theta)q^2
\end{array}\!\!\right)
+
\sin(\theta)
\left(\!\!\begin{array}{c}
\ttbeta(\theta)\\
h-\mu \cos(\theta)
\end{array}\!\!\right),
\end{equation}
on the cylinder $(\theta,q)\in S^1\times \R$, which is the same as $\{(m_3,r,q)\in\R^3: m_3^2+r^2=1\}$.

\paragraph{Wavetrains.} Steady states with vanishing $\xi$-derivative of $\theta$ and $q$ have $\varphi = q(x-st) +\Omega t$ and thus
correspond to the wavetrains discussed in \S\ref{s:wt} with wavenumber $k=q$ and frequency $\omega=sq - \Omega$. Hence, the ansatz \eqref{e:ansatz} removes all wavetrains whose wavenumber and frequency do not lie on the line $\{\omega= sk-\Omega\}$ in $(\omega,k)$-space. 

We may visualize this by drawing the line $\omega=sk-\Omega$ into the wavetrain existence and stability plots in the $(\omega,k)$-plane such as Figure~\ref{f:wt-stab}.  Specifically, for $s\neq 0$, steady states of \eqref{e:llg-coh} with $m_3\neq \pm1$ are wavetrains with wavenumber $q$ for which there exists $\theta$ such that 
\begin{equation}\label{e:qeq}
q= \frac{\Omega-\ttbeta(\theta)/\alpha}{s}. 
\end{equation}
In particular, for $\cc=0$ we have constant $\ttbeta(\theta)=\beta$ so that equilibria of \eqref{e:llg-coh} (other than $\pm\3$) have uniquely selected wavenumber $q$. Hence, heteroclinic solutions to \eqref{e:llg-coh} for $\cc=0$ can only be domain walls connecting $\pm\3$ or connect one of $\pm\3$ to a wavetrain. 

\paragraph{Coherent structure ODEs.} 
Writing \eqref{e:llg-coh} as an explicit ODE gives
\begin{equation}\label{e:ode}
\begin{array}{rl}
\theta' &= p\\
p' &= \sin(\theta) \left(h + (q^2-\mu) \cos(\theta) -(\Omega - sq) \right) - \alpha s p \\
q' &= \displaystyle  \alpha (\Omega-sq)  - \beta/(1+\cc\cos(\theta)) -\frac{s +2 \cos(\theta)q}{\sin(\theta)}p,
\end{array}
\end{equation}
whose study is the subject of the following sections. For later use we also note the `desingularization' by the (singular) coordinate change $p=\sin(\theta)\tp$ so that $\tp' = p'/\sin(\theta) - \tp^2 \cos(\theta)$, which gives
\begin{equation}\label{e:dode}
\begin{array}{rl}
\theta' &= \sin(\theta)\tp\\
\tp' &= h + (q^2-\mu) \cos(\theta) -(\Omega - sq) - \alpha s\tp - \cos(\theta)\tp^2 \\
q' &= \alpha (\Omega-sq)  - \beta/(1+\cc\cos(\theta)) - (s +2 \cos(\theta)q)\tp.
\end{array}
\end{equation}
Hence, \eqref{e:ode} is equivalent to \eqref{e:dode} except at zeros of $\sin(\theta)$. In particular, for $\tp=0$ the equilibria of \eqref{e:dode} with $\sin(\theta)\neq 0$ are those of \eqref{e:ode}, but $\theta= n\pi$, $n\in \mathbb{Z}$ are invariant subspaces which may contain equilibria with $\tp\neq 0$.

\bigskip
Next, we first consider various moving heteroclinic coherent structures with $s\neq 0$ and for $\cc=0$ give a complete analysis of stationary coherent structures ($s=0$). Coherent structures also emerge near the elliptic and hyperbolic wavetrain bifurcations that arise from fold points for $\cc\neq0$. However, the detailed analysis of this case is beyond the scope of this paper. 

\subsection{Homogeneous domain walls}

Classical domain walls connect antipodal equilibria at $x=\pm \infty$. For the model equation \eqref{eq:DWmotion} explicit (Walker) solutions are known to exist below a critical field $h$. These solutions exhibit a tilting of the azimuthal angle $\varphi=const.$ in order to balance precessional forces. An analogue situation arises in our context when $q=q'=0$. In this case we have $\varphi=\Omega t$ and therefore no spatially varying azimuthal profile.\footnote{After acceptance of the present manuscript for publication, we found that the sufficiency of $\cc=0$ for existence of such domains walls in Theorem~\ref{t:domwall} is contained in \cite{GRS:10}.} 

\begin{Theorem}\label{t:domwall}
Non-equilibrium coherent structure solutions with $q\equiv0$ for $s\neq 0$ and $|\cc|<1$ exist for $\mu<0$ and $\cc=0$ or $\beta=0$ only, and have $\Omega=\frac{h+\alpha\beta}{1+\alpha^2}$, $s^2 = -\frac{(\beta-\alpha h)^2}{\mu(1+\alpha^2)^2}$.  They are oscillating heteroclinic fronts connecting $\pm\3$ that solve $\theta'=\sigma\sqrt{-\mu}\sin(\theta)$ where $\sigma=\mathrm{sgn}(s(\alpha h-\beta))$. The family of such fronts is smooth and extends to $s=0$, where $h=\Omega=\beta/\alpha$, and fronts exist for both signs of $\sigma$. 
\end{Theorem}

\begin{Proof}
Recall $\ttbeta(\theta) = \beta/(1+\cc \cos(\theta))$ so $\frac{\rmd}{\rmd\xi}\ttbeta(\theta)=\cc \beta \sin(\theta)\theta'/(1+\cc\cos(\theta))^2$. 

Suppose a solution $(\theta,p,q)$ to \eqref{e:ode} has $q\equiv 0$, so also $q'\equiv 0$. Then the third equation of \eqref{e:ode} for $s\neq 0$ yields, using the first equation, $\theta' = \frac{\alpha\Omega-\ttbeta(\theta)}{s}\sin(\theta)$. Differentiation gives
\begin{align*}
\theta'' &=\frac{\alpha\Omega-\ttbeta(\theta)}{s}\cos(\theta)\theta' - \frac{\cc \beta \sin^2(\theta)}{s(1+\cc\cos(\theta))^2} \theta'\\ 
&= \sin(\theta)\frac{\alpha\Omega-\ttbeta(\theta)}{s^2}\left((\alpha\Omega-\ttbeta(\theta))\cos(\theta) - \frac{\cc \beta(1-\cos^2(\theta))}{(1+\cc\cos(\theta))^2})\right).
\end{align*}
On the other hand, the second equation of  \eqref{e:ode} requires
\[
\theta''=\sin(\theta)(h-\mu\cos(\theta)-\Omega-\alpha(\alpha\Omega-\ttbeta(\theta))).
\]
First consider $\cc=0$ or $\beta=0$ so that $\ttbeta(\theta)=\beta$ and these two right hand sides for $\theta''$ simplify. Equating them and comparing the coefficients of $\cos(\theta)^j$, $j=0,1$, yields $h=\Omega+\alpha(\alpha\Omega-\beta)$ and $\mu=- \left(\frac{\alpha\Omega-\beta}{s}\right)^2$, which means $(\alpha\Omega-\beta)/s = \sigma\sqrt{-\mu}$ for $\sigma=\sgn(s(\alpha\Omega-\beta))$. Taken together, the parameter conditions can be equivalently cast as the equations for $\Omega$, $s^2$ and $\sigma$ in the theorem statement. 

Hence, for $\cc=0$ or $\beta=0$ these parameter choices and $\theta' = \sigma\sqrt{-\mu}\sin(\theta)$ are necessary conditions for $q\equiv 0$. As a scalar equation, the only non-trivial and bounded solutions are heteroclinic orbits between equilibria. Taking $\Omega=\beta/\alpha + s\tilde\mu$ for some $\tilde\mu\neq 0$ gives a smooth parametrization up to $s=0$ and $\sigma=\sgn(\tilde\mu)$.

Conversely, for $\cc=0$ or $\beta=0$ and these choices of parameters, any $(\theta,p,q)(\xi)$ where $\theta(\xi)$ satisfies $\theta' = \sigma\sqrt{-\mu}\sin(\theta)$, $p=\theta'$ and $q\equiv0$ is a solution to \eqref{e:ode}.
 
\medskip
It remains to show that if $\beta\neq 0$ then $\cc=0$ is necessary for a non-trivial solution with $q\equiv 0$: Subtracting the two right hand sides for $\theta''$ from above and multiplication with $(1+\cc \cos(\theta))^3$ gives a polynomial in $\cos(\theta)$ of degree four. A straightforward computation shows that the 4th order term to vanish requires $\mu = - \alpha^2\Omega^2/s^2$. Using this the coefficients $a_j$ of $\cos(\theta)^j$, $j=0,1,2,3$, in this polynomial can be computed as
\begin{align*}
a_0  & = -h+(1+\alpha^2)\Omega -\alpha\beta +(\beta-\alpha\Omega)\beta \cc/s^2,\\
a_1 & = 3\cc((1+\alpha^2)\Omega-h)+\beta^2/s^2-\alpha\beta(2\cc+(2+\cc^2)\Omega/s^2),\\
a_2 & = \cc^2 (3 ((1 + \alpha^2) \Omega - h) - \alpha \beta ) - 3 \alpha \beta \Omega \cc/s^2,\\
a_3 & = \cc^3 ((1 + \alpha^2) \Omega - h) - \alpha \beta \Omega \cc^2/s^2.
\end{align*} 
We first solve $a_0=0$ trivially for $h$ and proceed with somewhat tedious, but straightforward calculations: substituting this $h$ into $a_1=0$ (which is linear in $\Omega$) we solve for $\Omega$, which uses $|\cc|<1$. Substituting the resulting $h,\Omega$ gives
\[
a_2=\frac{\cc\beta}{2s^2}\left(3\beta(\cc^2-1)+\alpha \cc s^2\right),
\]
so that for $a_2=0$ either $\cc=0$ (since $\beta\neq0$; note that then also $a_3=0$) or  $s^2=3\beta\frac{\cc^2-1}{\alpha\cc}$. In the latter case, substituting the previous $h,\Omega$ and this $s^2$ into $a_3$ would give $a_3 = \alpha\beta\cc^3/3\neq 0$. Hence, $\cc=0$ is necessary as claimed.
\end{Proof}


\begin{Remark}~
\begin{enumerate}
\item For $s=0$ further coherent structure solutions exist, but not as domain walls. See Theorem~\ref{t:coh-c0-q0} below. 
\item In \S\ref{s:fastfront} we find fast domain walls and fronts that have non-trivial $q$.
\item Numerical simulations suggest that these domain walls are dynamically stable solutions in the subsubcritical case. They are unstable in  the subcritical case $|h-\beta/\alpha|>-\mu$, which occurs for large $|s|$, since then either $\3$ or $-\3$ is unstable. 
\item The profile of these domain walls depends only on the parameter $\mu$. In particular, the subfamily parameterized by $h$ has arbitrarily large speed but constant shape, though the oscillation frequency $\Omega$ depends on $h$.
\end{enumerate}
\end{Remark}


\subsection{Moving front-type coherent structures}

Using the desingularized system \eqref{e:dode}, we prove existence of some non-stationary coherent structures of front-type, spatially connecting wavetrains or $\pm\3$.

\subsubsection{Near the fast limit $|s|\gg 1$}\label{s:fastfront}

\begin{Theorem}\label{t:coh-fast}
For any bounded set of $(\alpha, \beta, \mu, \Omega_0, \Omega_1)$ and $\cc\in(-1,1)$ there exists $s_0>0$ and neighborhood $U$ of $M_0:=\{\theta\in[0,\pi], \tp=q-\Omega_1=0\}$ such that for all $|s|\geq s_0$ the following holds for \eqref{e:dode} with $\Omega=\Omega_0 + \Omega_1 s$.  The heteroclinic orbits of \eqref{e:dode} in $U$ form a smooth family in the parameter $s^{-1}$ for each sign of $s$, which reverses their orientation. These heteroclinics and are in one-to-one correspondence with those of the ODE 
\begin{equation}\label{e:superslow}
\frac{\rmd}{\rmd\eta} \theta = -\frac{\alpha}{1+\alpha^2}\sin(\theta)\left(\frac{\ttbeta(\theta)}{\alpha}-h+(\Omega_1^2-\mu)\cos(\theta)\right),
\end{equation}
on the spatial scale $\eta = \xi/s$, which also gives the $\theta$-profile to leading order in $s^{-1}$. Moreover, for such a heteroclinic orbit $(\theta_{\rm h}, \tp_{\rm h}, q_{\rm h})(\xi)$ with $\theta^\sigma:=\lim_{\xi\to\sigma\infty} \theta_{\rm h}(\xi)\in\{0,\pi\}$ for $\sigma=1$ or $\sigma=-1$, the $q$-limit is/are 
\begin{equation}\label{e:frontq}
\lim_{\xi\to\sigma\infty} q_{\rm h}(\xi) = \Omega_1-\frac 1 s \left(\frac{h+\alpha\ttbeta(\theta^{\sigma})-\sigma\mu}{1+\alpha^2}-\Omega_0\right) + O(s^{-2}).
\end{equation}
In particular, for $\Omega_1\neq 0$ or $(1+\alpha^2)\Omega_0\neq h+\alpha\ttbeta(\theta^{\sigma})-\sigma\mu$ local wavenumbers are nontrivial: $q_{\rm h}\not\equiv 0$.
\end{Theorem}

Before proving the theorem we note the consequences of this for coherent structures and domains walls in in \eqref{e:ode} and \eqref{e:llg}.

\begin{Corollary}\label{c:fast}
The heteroclinic solutions of Theorem~\ref{t:coh-fast} are in one-to-one correspondence with heteroclinic solutions to \eqref{e:ode} and thus heteroclinic coherent structures in \eqref{e:llg} that lie in $U$ and connect $\theta=0,\pi$ or a wavetrain with $r\neq 0$. For $\theta\in(0,\pi)$ all properties carry over to \eqref{e:ode} with the bijection given by $p=\sin(\theta)\tp$. 
\end{Corollary}

\begin{Proof}
Recall that \eqref{e:dode} and \eqref{e:ode} are equivalent for $\theta\in(0,\pi)$. Since the limit of the vector field of \eqref{e:ode} along such a heteroclinic from \eqref{e:dode} is zero by construction in all cases. Hence, for each of the heteroclinic orbits in \eqref{e:dode} of Theorem~\ref{t:coh-fast}, there exist a heteroclinic orbit in \eqref{e:ode} in the sense of the corollary statement.
\end{Proof}

\begin{Corollary}
For any `bandgap' parameter set of \eqref{e:ode} such that there exist no wavetrains satisfying \eqref{e:qeq} for any $|s|>s_1$, for some $s_1>0$, there exist fast domain wall type coherent structures spatially connecting $\pm\3$ for all sufficiently large velocity $|s|$.
\end{Corollary}

\begin{Proof}
Choosing $\Omega_1=k$ there are by assumption no equilibria in \eqref{e:superslow} besides $\pm\3$, which are therefore connected by a heteroclinic orbit. Theorem~\ref{t:coh-fast} then implies the claim.
\end{Proof}

Such `bandgaps' occur in particular if $\mu>0$ for $\Omega_1^2\sim\mu$.

\begin{Remark}~
\begin{enumerate}
\item Concerning stability, Lemma~\ref{l:conststab} and Theorem~\ref{t:wt} imply that for $\cc=0$ the domain walls connecting $\pm\3$ might be stable in the subsubcritical case only since otherwise one of the asymptotic states is unstable: the unique wavetrain with $\theta\in(0,\pi)$ in the subsubcritical case and $\3$ or $-\3$ in the sub- and supercritical cases. However, it may be that some fronts are stable in a suitable weighted sense as invasion fronts into an unstable state. 
\item For increasing speeds these solutions are decreasingly localized, hence far from a sharp transition.  
\item The uniqueness statement in the corollaries is limited, since in the $(\theta,p,q)$-coordinates the neighborhood $U$ from the theorem is `pinched' near $\theta=0,\pi$: a uniform neighborhood in $(\theta,\tp)$ has a sinus-shaped boundary in $(\theta,p)$. 
\end{enumerate}
\end{Remark}

\begin{Remark}
Part of the family homogeneous domains walls from Theorem~\ref{t:domwall}, where $\cc=0$, is a continuation to smaller $|s|$ of homogeneous ($q\equiv0$) fronts in the family of Theorem~\ref{t:coh-fast}. The latter are decreasingly localized, which requires in the former that $\sqrt{\mu}=\calO(s)$. Specifically, $\mu=-\frac{(\beta-\alpha h)^2}{s^2(1+\alpha^2)^2}$ and $\Omega_0=\frac{h+\alpha\beta}{1+\alpha^2}$, $\Omega_1=0$ in the heteroclinics of Theorem~\ref{t:coh-fast}. Then $\mu\to 0$ as $s^2\to\infty$ so that $\mu=0$ in the leading order equation \eqref{e:superslow} and in \eqref{e:frontq} $\mu$ is removed from the order $s^{-1}$ term. Since  $\sigma s\sqrt{-\mu} = -(\beta-\alpha \Omega_0)$ and $\beta-\alpha\Omega = -\frac{\alpha}{1+\alpha^2}\left(h-\frac{\beta}{\alpha}\right)$ indeed  \eqref{e:superslow} equals the equation in Theorem~\ref{t:domwall}. In particular,  \eqref{e:frontq} is consistent with $q\equiv0$.
\end{Remark}


Finally, remark that the ODE \eqref{e:superslow} is the spatial variant of the temporal heteroclinic connection in \eqref{e:llg-cyl}: setting all space derivatives to zero, the $\theta$-equation of \eqref{e:llg-cyl} reads
\[
-\partial_t \theta = \frac{\alpha}{1+\alpha^2}\sin(\theta)\left(h-\frac{\ttbeta(\theta)}{\alpha}-\mu\cos(\theta)\right),
\]
which is \eqref{e:superslow} with $\mu$ replaced by $\Omega_1^2-\mu$ and up to possible direction reversal. 
Since $\Omega_1=q$ on the slow manifold $M_0$ (i.e. at leading order), the reduced flow equilibria reproduce the wavetrain existence condition \eqref{e:psi}.
This kind of relation between temporal dynamics and fast travelling waves holds formally (but in general not rigorously)  for any evolution equation in one space dimension. Here the symmetry makes the temporal ODE scalar.

\begin{Proof}[ (Theorem~\ref{t:coh-fast})]
Let us set $s=\eps^{-1}$ so that the limit to consider is $\eps\to 0$. Since we will rescale space with $\eps$ and $\eps^{-1}$ this means sign changes of $s$ reflect the directionality of solutions.

The existence proof relies on a geometric singular perturbation argument and we shall use the terminology from this theory, see \cite{Fenichel,Jones}, and also sometimes suppress the $\eps$-dependence of $\theta,\tp,q$.

Upon multiplying the $\tp$- and $q$-equations of \eqref{e:dode} by $\eps$ we obtain the, for $\eps\neq 0$ equivalent, `slow' system 
\begin{equation}\label{e:s-dode}
\begin{aligned}
\theta' &= \sin(\theta)\tp\\
\eps\tp' &= -\alpha \tp + q - \Omega_1 + \eps(h + (q^2-\mu) \cos(\theta) -\Omega_0 - \cos(\theta)\tp^2) \\
\eps q' &= -\tp - \alpha (q- \Omega_1) + \eps(\alpha\Omega_0  - \ttbeta(\theta) -2 \cos(\theta)q \tp).
\end{aligned}
\end{equation}

Setting $\eps=0$ gives the algebraic equations 
\[
A\begin{pmatrix}\tp\\q\end{pmatrix}=-\Omega_1\begin{pmatrix}-1\\\alpha\end{pmatrix}, \mbox{ where }
A =- 
\begin{pmatrix}
\alpha & -1\\
1 & \alpha 
\end{pmatrix}.
\]
Since $\det A = 1+\alpha^2>0$ the unique solution is $\tp=q-\Omega_1=0$ and thus the `slow manifold' is $M_0$ as defined in the theorem, with `slow flow' given by
\[
\theta' = \sin(\theta)\tp.
\]
Since $\tp = 0$ at $\eps=0$, $M_0$ is a manifold (a curve) of equilibria at $\eps=0$, so that the slow flow is in fact `superslow' and will be considered explicitly below. Since the slow manifold is one-dimensional (and persists for $\eps>0$ as shown below) it suffices to consider equilibria for $\eps>0$. These lie on the one hand at $\theta=\theta_0\in\{0,\pi\}$, if  
\[
A\begin{pmatrix}\tp\\q\end{pmatrix} + \Omega_1\begin{pmatrix}-1\\\alpha\end{pmatrix} + \eps F(\tp,q)=0\;,\quad F(\tp,q) := 
\begin{pmatrix}
h + \sigma(q^2-\mu)  -\Omega_0 - \sigma\tp^2 \\
\alpha\Omega_0  - \ttbeta(\theta_0) -2 \sigma q\tp
\end{pmatrix},
\]
where $\sigma = \cos(\theta_0)\in \{-1,1\}$. Since $\det A = -(1+\alpha^2)<0$ the implicit function theorem provides a locally unique curve of equilibria $(\tp_\eps,q_\eps)$ for sufficiently small $\eps$, where 
\[
\left.\frac{\rmd}{\rmd\eps}\right|_{\eps=0} 
\begin{pmatrix}
\tp_\eps\\
q_\eps
\end{pmatrix}
 = -A^{-1}F(0,0) = -A^{-1}
\begin{pmatrix}
h-\sigma\mu-\Omega_0\\
\alpha\Omega_0  - \beta^\pm
\end{pmatrix}.
\]
This proves the claimed location of asymptotic states.

On the other hand, for $\theta\neq 0$ system \eqref{e:dode} is equivalent to \eqref{e:ode}. From the previous considerations of equilibria (=wavetrains) we infer that the unique equilibria in an $\eps$-neighborhood of $M_0$ are those at $\theta=\theta_0$, $(\tp,q) = (\tp_\eps,q_\eps)$ together with the 
possible additional $\theta$-values of wavetrains, where $k$ is now replaced by $q=\eps(\Omega_0-\ttbeta(\theta)/\alpha)$.

\medskip
Towards the persistence of $M_0$ as a perturbed invariant manifold for $|\eps|>0$, let us switch to the `fast' system by rescaling the time-like variable to $\zeta = \xi/\eps$. With $\dot \theta = \rmd \theta /\rmd \zeta$ etc., this gives
\begin{equation}\label{e:f-dode}
\begin{aligned}
\dot\theta &= \eps\sin(\theta)\tp\\
\dot\tp &= -\alpha \tp + q-\Omega_1 + \eps(h + (q^2-\mu) \cos(\theta) -\Omega_0 - \cos(\theta)\tp^2) \\
\dot q &= -\tp - \alpha (q-\Omega_1) + \eps(\alpha\Omega_0  - \ttbeta(\theta) -2 \cos(\theta)q \tp).
\end{aligned}
\end{equation}
Note that $M_0$ is (also) a manifold of equilibria at $\eps=0$ in this system and the linearization of \eqref{e:f-dode} in $M_0$ for transverse directions to $M_0$ is given by $A$. Since the eigenvalues of $A$, $-\alpha\pm\rmi$, are away from the imaginary axis, $M_0$ is normally hyperbolic and therefore persists as an $\eps$-close invariant one-dimensional manifold $M_\eps$, smooth in $\eps$ and unique in a neighborhood of $M_0$. See \cite{Fenichel}. The aforementioned at least two and at most three equilibria lie in $M_\eps$, and, $M_\eps$ being one-dimensional, these must be connected by heteroclinic orbits. 

\medskip
For the connectivity details it is convenient to derive an explicit expression of the leading order flow. We thus switch to the superslow time scale $\eta = \eps\xi$ and set $\op = \tp/\eps$, $\oq=(q-\Omega_1)/\eps$, which changes \eqref{e:s-dode} to (subdindex $\eta$ means $\rmd/\rmd\eta$)
\begin{equation}\label{e:ss-dode}
\begin{aligned}
\theta_\eta &= \sin(\theta)\op\\
\eps\op_\eta &= -\alpha \op + \oq + h +(\Omega_1^2-\mu)\cos(\theta) -\Omega_0 + \eps\cos(\theta)(\eps (\oq^2- \op^2) +2\oq\Omega_1) \\
\eps \oq_\eta &= -\op - \alpha \oq + \alpha\Omega_0  - \ttbeta(\theta) -2\eps \cos(\theta)(\eps \oq \op + \Omega_1\op).
\end{aligned}
\end{equation}
At $\eps=0$, solving the algebraic equations for $(\op,\oq)$ gives
\[
\begin{pmatrix}
\op\\
\oq
\end{pmatrix}
=-A^{-1} 
\begin{pmatrix}
h+(\Omega_1^2-\mu)\cos(\theta)-\Omega_0\\
\alpha\Omega_0-\ttbeta(\theta)
\end{pmatrix}
= \frac{1}{1+\alpha^2}
\begin{pmatrix}
\alpha h+ \alpha(\Omega_1^2-\mu)\cos(\theta) - \ttbeta(\theta)\\
(1+\alpha^2)\Omega_0-h - (\Omega_1^2-\mu)\cos(\theta) - \alpha\ttbeta(\theta)
\end{pmatrix}
\]
whose first component gives $\op$ so that the leading order superslow flow on the invariant manifold is indeed given by \eqref{e:superslow}.
\end{Proof}

\subsubsection{The case of small amplitudes}

In this section we consider small amplitude coherent structures, which means $q$ must lie near a bifurcation point of wavetrains. Here we focus on the intersection points of the solution curves from \eqref{e:wt-om} with $\theta=\theta_0=0,\pi$, which gives 
\begin{equation}\label{e:inters}
\cos(\theta_0) \left(q^2-\mu\right) = \frac{\ttbeta(\theta_0)}{\alpha} - h.
\end{equation}
For $\cc=0$ this is possible for super- and subcritical anisotropy only, compare Figure~\ref{f:wavetrains}.

We show that these intersection points are pitchfork-type bifurcations in \eqref{e:dode} that give rise to front-type coherent structures. As in the previous section, we locate such solutions in \eqref{e:ode} from an analysis of \eqref{e:dode}. 

It is convenient to write \eqref{e:wt-om} in terms of $m_3=\cos(\theta)$ so equilibria of \eqref{e:dode} with $\tp=0$ solve 
\[
\tGamma(m_3):=\frac{\tbeta(m_3)}{\alpha} - \left(q^2(m_3)-\mu\right)m_3 - h=0,
\]
where $q(m_3)$ is the selected $q$ from \eqref{e:qeq}. Recall $\tbeta(m_3) = \frac{\beta}{1+\cc m_3}$.

\begin{Theorem}\label{t:coh-small}
Consider $\theta=\theta_0\in\{0,\pi\}$ and set $\mz:=\cos(\theta_0)$. Suppose that parameters of \eqref{e:dode} are such that $s\neq 0$, $\tGamma(\mz)=0$ and $\partial_{m_3}\tGamma(\mz) 
\neq 0$.  Then the equilibrium point $(\theta_0, 0, q(\mz))$ of \eqref{e:dode} undergoes a pitchfork bifurcation upon any perturbation of $h$ or $\mu$. 

More precisely, let $S_\eps=(\alpha_\eps,\beta_\eps, h_\eps,\mu_\eps,\Omega_\eps,\cc(\eps),s_\eps)$,  $\eps\in(-\eps_0,\eps_0)$ for some $\eps_0>0$, be a curve in the parameter space of \eqref{e:dode} with $|\cc(\eps)|<1$, $\alpha_\eps>0$, $s_\eps\neq 0$ and such that $S_0$ satisfies $\tGamma(\mz)=0$ and $\tgamma:=\partial_{m_3}\tGamma(\mz)\neq0$.

Then \eqref{e:dode} has a curve of equilibria  $(\theta_0,\tp_\eps,q_\eps)$, with possibly adjusted $\eps_0$, such that $\tp_0=0, q_0=q(\mz)$ and two equilibria with $\theta\neq \theta_0$ bifurcate from $(\theta_0, 0, q(\mz))$ for increasing $\eps$ if, with parameters $S_\eps$, $\mz\tgamma\partial_\eps\tGamma(\mz)|_{\eps=0}>0$. The bifurcating equilibria are connected to $(\theta_0,\tp_\eps,q_\eps)$ by heteroclinic orbits which converge to $(\theta_0,\tp_\eps,q_\eps)$ as $\ell \xi\to\infty$ for $\ell=-\mz \tp_\eps$. 

Specifically, this occurs if $h_\eps=h_0-\mz\tgamma\eps$ or $\mu_\eps=\mu_0+\tgamma \eps$ and $\ell=\tgamma s$, or if $\beta_\eps = \beta_0 + \mz\tgamma \eps$ and $\ell = -(2q(\mz)\mz+s)\tgamma$, with all other parameters fixed in each case.
\end{Theorem}

Analogously to Corollary~\ref{c:fast} we have

\begin{Corollary}\label{c:coh-small}
The heteroclinic solutions of Theorem~\ref{t:coh-small} are in one-to-one correspondence with heteroclinic solutions to \eqref{e:ode}, connecting to $\theta\equiv0$ or $\theta=\pi$ in $U$. Bounded solutions for $\theta\not\in\{0,\pi\}$ are also in one-to-one correspondence.\end{Corollary}

\begin{Proof}[ (Theorem~\ref{t:coh-small})]
Note that $\tGamma(m_3)=0$ with $|m_3|<1$ is equivalent to (and if $|m_3|=1$ sufficient for) the existence of an equilibrium of \eqref{e:dode} at $\theta$ with $\cos(\theta)=m_3$, $q=q(m_3)$ from \eqref{e:inters} and $\tp=0$. Assuming $\tgamma=\partial_{m_3}\tGamma(\mz)\neq0$ and $\partial_\nu \tGamma(\mz)\neq0$ for $\nu=h$ or $\nu=\mu$ implies existence of a locally unique curve of equilibria $m_3(\nu)$ that transversely crosses $\mz$. 
The case of parameters $S_\eps$ is analogous with a curve $m_3(\eps)$, where $\partial_\eps m_3(0) = -\partial_\eps \tGamma(\mz)/\partial_{m_3}\tGamma(\mz)|_{\eps=0}$ having the sign of $-\mz$ means bifurcation of two equilibria for $\eps>0$.

It remains to show that the center manifold associated to the bifurcation is one-dimensional, and to obtain the directionality of heteroclinics. 

For the former it suffices to show that the linearization at the bifurcation point has only a simple eigenvalue on the imaginary axis, namely at zero. The linearization of \eqref{e:dode} in any point with $\tp=0$, $\theta=\theta_0$ gives the $3\times 3$ matrix
\[
\tilde A= 
\left(\begin{array}{c|c}
0  & 0\; 0\\
\hline
\begin{array}{c}0\\0\end{array} & B
\end{array}\right), \;
B = 
\left(\begin{array}{cc}

-\alpha s  & 2q \mz+s\\
-(s+2q\mz) & -\alpha s
\end{array}\right),
\] 
which has a kernel with eigenvector $(1,0,0)^t$. The remaining eigenvalues are those of $B$, which are $-\alpha s\pm(s+\mz 2q)\rmi$. Since these lie off the imaginary axis for $s\neq0$ there is indeed at most one simple zero eigenvalue on the imaginary axis. This implies the existence of a one-dimensional center manifold which includes all equilibria and heteroclinic connections near $(\theta_0,0,q(\mz))$ for nearby parameters. 
Equilibria in the symmetry plane $\{m_3=\mz\}$ are solutions of $(\tp',q')=:G(\tp,q)=0$ with $G$ given by \eqref{e:dode}. Since $DG(0,q(\mz))=B$ is invertible we obtain a curve $(\theta_0,\tp_\eps,q_\eps)$ for parameters at $S_\eps$ as claimed.

The uniqueness of bifurcating equilibria on either side of the symmetry plane and invariance of the one-dimensional center manifold implies existence and local uniqueness of heteroclinic connections for $\sgn(\eps) = \sgn(\mz\tgamma\partial_\eps\tGamma(\mz)|_{\eps=0})$. In order to determine the directionality of these, a perturbation in the kernel gives $\theta'=\sin(\theta_0+\delta)\tp = \mz\delta\tp + \calO(\delta^2)$ so that for $\mz\tp_\eps<0$ the equilibrium at $(\theta_0,\tp_\eps,q_\eps)$ is stable in the center manifold, and unstable for reversed sign. Since existence of heteroclinics requires $\sgn(\eps) = \sgn(\mz\tgamma\partial_\eps\tGamma(\mz)|_{\eps=0})$ this implies stability if $\ell:=-\tgamma\partial_\eps\tGamma(\mz)\partial_\eps\tp|_{\eps=0}<0$ and thus convergence to $(\theta_0,\tp_\eps,q_\eps)$ as $\ell\xi\to \infty$.

For $h_\eps=h_0-\mz\tgamma\eps$ with otherwise fixed parameters $\partial_\eps \tGamma(\mz) = \mz\tgamma$ so heteroclinics exist for $\eps>0$. On the other hand, $\partial_\eps \tp_\eps|_{\eps=0}$ is the first component $-\mz\tgamma\alpha s \det(B)$ of $-B^{-1} \partial_h G(0,q(\mz))(-\mz\tgamma)$, where $\alpha\det(B)>0$. Hence, $\ell = s\tgamma$ as claimed. The cases $\eps=\mu, \beta$ are determined analogously using $\partial_\mu \tGamma(\mz) = \mz$, $\partial_\mu G(0,q(\mz) = (-\mz,0)$ and $\partial_\beta \tGamma(\mz) = ((1+\cc\mz)\alpha)^{-1}>0$, $\partial_\beta G(0,q(\mz) = (0,-((1+\cc\mz)\alpha)^{-1})$.
\end{Proof}

\subsection{Stationary coherent structures for $\cc=0$}\label{s:stat-coh}

In this section we consider the case $s=0$ (which does not imply time-independence) and $\cc=0$ (which will imply integrability), so that equations \eqref{e:ode} reduce to
\begin{equation}\label{e:ode-c0}
\begin{array}{rl}
\theta'' &= \sin(\theta) \left(h-\Omega + (q^2-\mu) \cos(\theta)\right) \\
q' &=  \alpha \Omega- \beta -2\cot(\theta)\theta'q.
\end{array}
\end{equation}

In case $\Omega\neq\beta/\alpha$ there are no equilibria and it will be shown at the end of this section that there are no coherent structure-type solutions in that case. Recall from \S\ref{s:symc0} that for $\cc=0$ we may choose coordinates so that $\beta=0$, which means $\Omega=0$ and thus stationary coherent structures are turned into standing waves. However, we choose not to remove the parameter $\beta$ in order to emphasize the typically oscillatory nature of solutions to \eqref{e:llg} and for consistency in parameter relations. Nevertheless, the symmetries and integrals that we will find are  consequences of this reducibility.

For $\Omega=\beta/\alpha$, system \eqref{e:ode-c0} is invariant under the reflection $q\to -q$ so that $\{q=0\}$ is an invariant plane which separates the three dimensional phase space. In particular, \emph{there cannot be connections between equilibria (=wavetrains) with opposite signs of $q$, that is, sign reversed spatial wavenumbers}. 

\subsubsection{Homogeneous solutions ($q=0$)} 

\begin{figure}
\begin{center}
\includegraphics[width=0.24\textwidth]{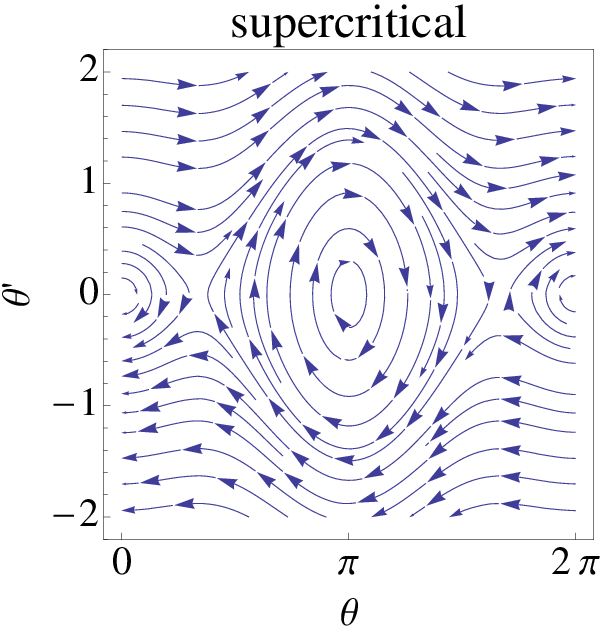}
\includegraphics[width=0.24\textwidth]{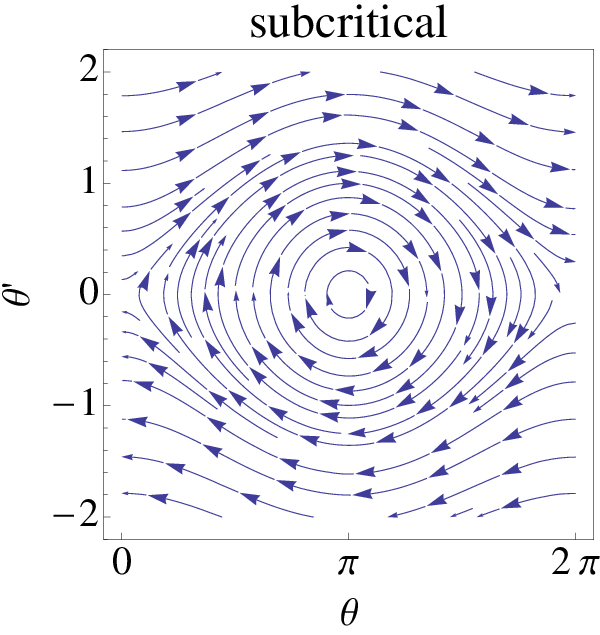}
\includegraphics[width=0.24\textwidth]{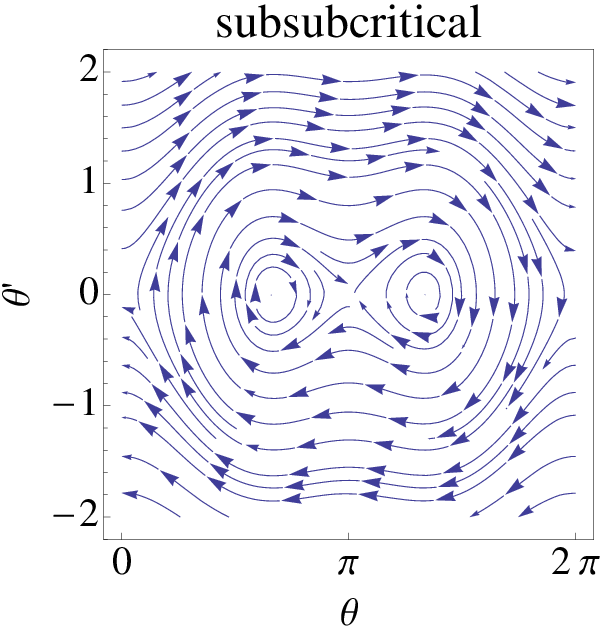}
\includegraphics[width=0.24\textwidth]{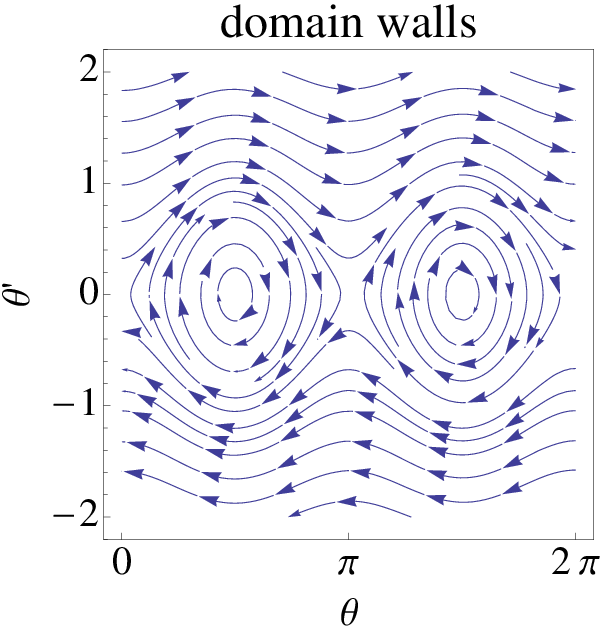}\\
\hfill (a)\hfill~\hfill (b)\hfill~\hfill (c)\hfill~\hfill (d)\hfill~
\caption{Phase plane streamplots of \eqref{e:ode-c0-q0} with \textsc{Mathematica}. (a)-(c) have $h-\Omega=1/2$. (a) supercritical anisotropy (here $\mu=1$), (b) subcritical (here $\mu=0$), (c) subsubcritical (here $\mu=-1$), (d) subsubcritical case that allows for standing domain walls, $h=\Omega$, $\mu=-1$.}
\label{f:phaseplane-c0-q0}
\end{center}
\end{figure}

Solutions in the invariant set $\{q=0\}$ have the form $m(\xi) = r(\xi)\exp(\rmi t \Omega )$ and \eqref{e:llg-coh} turns into a second order ODE on the circle $\{m_3^2+r^2=1\}$. The ODE for $\theta$ from \eqref{e:ode-c0} is given by the nonlinear pendulum equation
\begin{equation}\label{e:ode-c0-q0}
\theta'' = \sin(\theta) \left(h-\Omega-\mu \cos(\theta)\right),
\end{equation}
which is invariant under $\theta\to-\theta$ and is Hamiltonian with potential energy 
\[
P_0(\theta)=\cos(\theta)(h-\Omega- \frac\mu 2 \cos(\theta)).
\]
The symmetry \eqref{e:parasym} applies and we therefore assume in the following that $\Omega=\beta/\alpha <h$. 

We plot the qualitatively different vector fields of \eqref{e:ode-c0-q0} in Figure~\ref{f:phaseplane-c0-q0} and some profiles in Figure~\ref{f:profiles-c0-q0}. Coherent structure solutions are completely characterized via the figure, which we formulate next explicitly for the original PDE with $\m=(m,m_3)$, $m=r \rme^{\rmi\varphi}$, $r=\sin(\theta)$, $m_3 = \cos(\theta)$. Homoclinic profiles may be interpreted as (dissipative) solitons. The heteroclinic connections in item 2(a) can be viewed as (dissipative) solitons with `phase slip'.

\begin{Theorem}\label{t:coh-c0-q0} Let $s=0$ and $\Omega=\beta/\alpha$ and consider solutions to \eqref{e:llg-mod} of the form \eqref{e:ansatz} with $\varphi$ constant in $\xi$, i.e., $q=0$. These oscillate in time pointwise about the $\3$-axis with frequency $\Omega=\beta/\alpha$. Assume without loss of generality, due to \eqref{e:parasym}, that $h>\Omega$.

\begin{enumerate}
\item Subcritical anisotropy $h-\Omega >|\mu|>0$. There exist no nontrivial wavetrains with $q=0$, and the coherent structure solutions with $q=0$ are a pair of homoclinic profiles to $\hat e_3$, and three one-parameter families of periodic profiles, one bounded and two semi-unbounded. The homoclinic profiles each cross once through $-\3$ in opposite $\theta$-directions. The limit points of the bounded curve of periodic profiles are $-\3$ and the union of homoclinic profiles. Each of the homoclinics is the limit point of one of the semi-unbounded families, each of which has unbounded $\theta$-derivatives. The profiles from the bounded family each cross $-\3$ once during a half-period, the profiles of the unbounded family cross both $\pm\3$ during one half-period.
\item Suppose super- or subsubcritical anisotropy $|\mu|>h-\Omega$. There exists a wavetrain with $k=0$, which is stable in the supercritical case ($\mu>0$) and unstable in the subsubcritical case ($\mu<0$). In \eqref{e:ode-c0-q0} this appears in the form of two equilibria being the symmetric pair of intersection points of the wavetrain orbit and a meridian on the sphere, phase shifted by $\pi$ in $\varphi$-direction. Details of the following can be read off Figure~\ref{f:phaseplane-c0-q0} analogous to item 1.

\begin{enumerate}
\item In the supercritical case ($\mu>h-\Omega$) the coherent structure solutions with $q=0$  are two pairs of heteroclinic connections between the wavetrain and its phase shift, and four curves of periodic profiles; two bounded and two semi-unbounded. 
\item In the subsubcritical case ($\mu<h-\Omega$) the coherent structure solutions with $q=0$  are two pairs of homoclinic connections to $\pm\3$, respectively, and five curves of periodic profiles, three bounded and two semi-unbounded. 
\end{enumerate}

\item The degenerate case $h=\Omega$, $\mu<0$ is the only possibility for profiles of stationary coherent structures to connect between $\pm\3$, which then come in a pair as in the corresponding panel of Figure~\ref{f:phaseplane-c0-q0}. The remaining coherent structures with $q=0$ are analogous to the supercritical case with $\pm\3$ and the pair of wavetrain and its phase shift interchanged.
\end{enumerate}
\end{Theorem} 

\begin{Proof}
As for wavetrains discussed in \S\ref{s:wt}, the condition $\cos(\theta) = (h-\Omega)/\mu$ yields the existence criterion $|h-\Omega|<|\mu|$ for an equilibrium to \eqref{e:ode-c0-q0} in $(0,\pi)$. The derivative of the right hand side of \eqref{e:ode-c0-q0} at $\theta=0$ is $h-\Omega-\mu$, which dictates the type of all equilibria and only saddles generate heteroclinic or homoclinic solutions.

It remains to study the connectivity of stable and unstable manifolds of saddles, which is given by the difference in potential energy $P_0(\theta)$. Since $P_0(0)-P_0(\pi) = 2(h-\Omega)$ the claims follow.
\end{Proof}

\subsubsection{Non-homogeneous solutions ($q\neq 0$)}

\begin{figure}
\begin{center}
\begin{tabular}{cc}
\begin{minipage}{0.45\textwidth}
\includegraphics[width=0.8\textwidth]{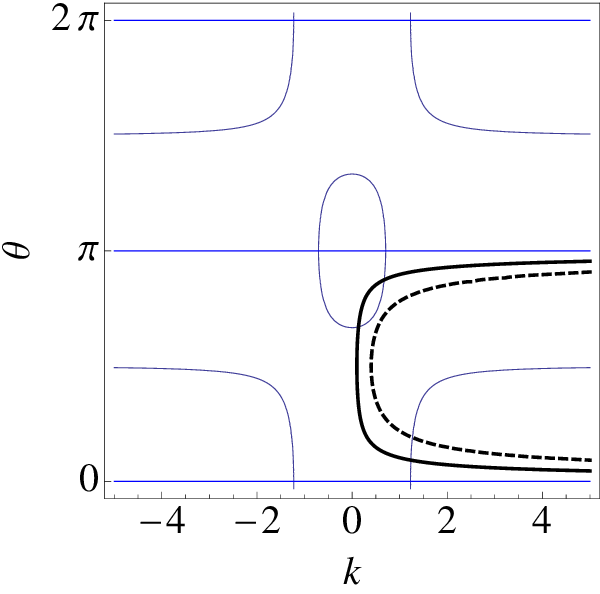}
\end{minipage}
&
\begin{minipage}{0.4\textwidth}
\begin{center}
\includegraphics[width=0.7\textwidth]{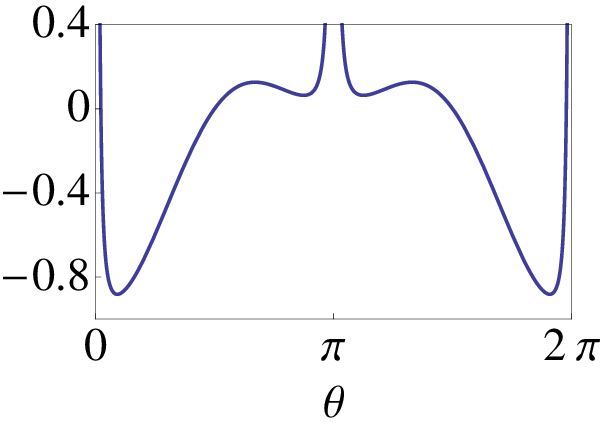}\\
\includegraphics[width=0.7\textwidth]{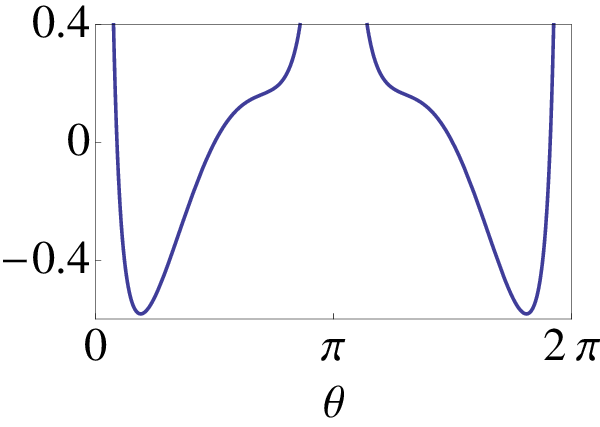}
\end{center}
\end{minipage}\\
(a) & (b)
\end{tabular}
\caption{(a) Figure~\ref{f:wavetrainstab}(a) with solutions to \eqref{e:coh-c0-q} for $\theta\in(0,\pi)$ and $C=0.1$ (thick solid line), and $C=0.4$ (thick dashed line). (b) Upper panel: potential $P(\theta)$ for parameters as in (a) and $C=0.1$. Lower panel: same with $C=0.4$.}
\label{f:wt-cyl-q}
\end{center}
\end{figure}

In order to study \eqref{e:ode-c0} for $q\neq 0$, we note the following first integral. Since $\Omega=\beta/\alpha$, the equation for $q$ can be written as 
\[
(\log|q|)' = -2(\log|\sin(\theta)|)',
\]
and therefore explicitly integrated. With integration constant $C = \sin(\theta(0))^2 |q(0)|$ this gives 
\begin{equation}\label{e:coh-c0-q}
q = \frac{C}{\sin(\theta)^2}.
\end{equation}
Substituting this into the equation for $\theta$ yields the nonlinear pendulum
\begin{equation}\label{e:ode-c0-red}
\theta'' = \sin(\theta) \left(h-\Omega-\mu \cos(\theta)\right) + C^2\frac{\cos(\theta)}{\sin(\theta)^3},
\end{equation}
with singular potential energy
\[
P(\theta) = P_0(\theta) + \frac 1 2 C^2 \cot(\theta)^2.
\]
The energy introduces barriers at multiples of $\pi$ so that solutions for $q\neq 0$ cannot pass $\pm\hat e_3$, and as $C$ increases the energy landscape becomes qualitatively independent of $\Omega, h, \mu$. 

As an immediate consequence of the energy barriers and the fact that only stable wavetrains are saddle points we get

\begin{Theorem}\label{t:coh-c0}
Let $s=0$ and $\Omega=\beta/\alpha$. Consider solutions to \eqref{e:llg-mod} of the form \eqref{e:ansatz}. Intersections in $(q,\theta)$-space of the curve $\mathcal{C}$ given by \eqref{e:coh-c0-q} with the wavetrain existence curves $\mathcal{W}$ from \eqref{e:psi} (with $k$ replaced by $q$) are in one-to-one correspondence with equilibria of \eqref{e:ode-c0}. 

Consider supercritical anisotropy $0\leq h-\Omega < \mu$ and assume $C$ is such that $\mathcal{C}$ transversely intersects the component of $\mathcal{W}$ which intersects $\{q=0\}$. Then the intersection point with smaller $q$-value corresponds to a spectrally stable wavetrain, and in \eqref{e:ode-c0} there is a pair of homoclinic solutions to this wavetrain. All other intersection points are unstable wavetrains.

All other non-equilibrium solutions of the form \eqref{e:ansatz} with $s=0$ are periodic in $\xi$, and this is also the case for all other parameter settings.
\end{Theorem}

The homoclinic orbit is a soliton-type solution to \eqref{e:llg} with asymptotic state a wavetrain (cf.\ Figure~\ref{f:cohex}).

Notably, the tangential intersection of $\mathcal{C}$ and $\mathcal{W}$ is at the sideband instability.

In Figure~\ref{f:wt-cyl-q} we plot an illustration in case of supercritical anisotropy $\mu>\Omega-h>0$. In the upper panel of (b) the local maxima each generate a pair of homoclinic solutions to the stable wavetrain it represents. The values of $q$ that it visits lie on the bold curve in (a), whose intersections with the curve of equilibria are the local maxima and minima. The lower panel in (b) has $C=0.4$ and the local maxima disappeared. All solutions are periodic and lie on the thick dashed curve in (a).

\subsubsection{Absence of wavetrains ($\Omega\neq \beta/\alpha$)}

In this case the first integral $Q= \log(|q|\sin(\theta)^2)$ is monotone,  
\[
Q' = \frac{\Omega-\beta/\alpha}{q} \neq 0,
\]
and therefore $|q|$ is unbounded as $\xi\to\infty$ or $\xi\to-\infty$. In view of \eqref{e:ode-c0-red}, we also infer that $\theta$ oscillates so that there are no relevant solutions of the form \eqref{e:ansatz}.

\bibliographystyle{acm}
\bibliography{wavetrain}

\end{document}